\documentclass[11pt]{article}
\usepackage[T1]{fontenc}
\usepackage[utf8]{inputenc}
\usepackage[textwidth=6in,textheight=8in]{geometry}
\usepackage{amsmath,amssymb,amsthm,mathtools,booktabs,array,graphicx}
\usepackage{microtype}
\usepackage[numbers,sort&compress]{natbib}
\usepackage{xurl}
\usepackage[hidelinks]{hyperref}
\usepackage{placeins}

\newcommand{\R}{\mathbb{R}}
\newcommand{\Nzero}{\mathbb{N}_0}
\newcommand{\Xdot}{\dot{\mathbf X}}
\newcommand{\Thetahat}{\widehat{\mathbf\Theta}}
\DeclareMathOperator{\diag}{diag}
\newtheorem{proposition}{Proposition}[section]
\newtheorem{theorem}[proposition]{Theorem}
\newtheorem{remark}[proposition]{Remark}

\title{Data-Driven Discovery and Rigorous Certification of Periodic Solutions in Infinite-Dimensional Dynamical Systems}
\author{Aaron Loh\thanks{American Heritage Schools, Palm Beach, FL, USA. Email: \href{mailto:aaronsloh@gmail.com}{aaronsloh@gmail.com}}
\and J. D. Mireles James\thanks{Florida Atlantic University, Boca Raton, FL, USA. Email: \href{mailto:jmirelesjames@fau.edu}{jmirelesjames@fau.edu}}}
\date{November 2026}

\begin{document}
\maketitle

\begin{abstract}

In this work, we consider a standard regression technique for
time series analysis which fits a low dimensional proxy system
to a higher (even infinite) dimensional target model. We then
ask how well the long time dynamics of the proxy describe those
of the target. Standard short-term tracking diagnostics cannot
resolve this question, and we propose a novel diagnostic that
exploits a-posteriori computer-assisted proofs. That is, we
locate periodic solutions in the proxy model, and use them to
generate approximate solutions for the target. We then try to
certify the existence of a true periodic solution of the
target, nearby. When successful, this certification procedure
provides mathematically rigorous information about how well
an invariant structure in the proxy approximates a corresponding
invariant structure in the target model. We show that this
procedure can succeed without the application of differential
corrections in the target system.

To illustrate the utility of the procedure, we develop
surrogates using time series data from the Mackey--Glass delay
differential equation, and for the Kuramoto--Sivashinsky (KS)
partial differential equation, using periodic boundary
conditions for forecasting and Neumann boundary conditions
for certification.
When considered as dynamical systems, each of these models has
an infinite dimensional phase space, and each is well known to
exhibit rich dynamics. In both test cases we are able to
extract invariant dynamics from low dimensional proxies and
use this to validate the existence of corresponding structures
in the targets. For KS we also learn proxies from fixed-sensor
time series; their candidates converge to true periodic
solutions, one certified after Newton refinement, but the
regression becomes unstable before an uncorrected proof is reached.
These results demonstrate assessment ability
beyond short-term forecasting.

\end{abstract}

\section{Introduction}
\label{sec:introduction}

Data-driven dimension reduction is essential for analyzing complex systems, and tools 
like principal component analysis \cite{jolliffe2016}, multidimensional scaling 
\cite{borg2005}, and nonlinear manifold-learning techniques \cite{lee2007} provide 
low-dimensional representations of high-dimensional data in diverse applications 
throughout the natural and social sciences. For time-dependent data sets, it is 
natural to consider dimension reduction techniques rooted in dynamical systems theory, 
which aim to capture not only the geometry of the data but also describe its temporal 
evolution. Modern approaches such as Koopman operator analysis \cite{brunton2022koopman}, 
nonlinear model reduction using spectral submanifolds \cite{haller2025}, and operator 
inference for reduced-order dynamical systems \cite{peherstorfer2016} illustrate this 
principle. The present work considers questions raised by the data-driven vector-field 
reconstruction techniques like those of Crutchfield--McNamara and Cremers--H\"ubler \cite{cremers1987,crutchfield1987}.

To be precise, suppose that $f \colon \mathbb{R}^n \to \mathbb{R}^n$ is a low-dimensional vector 
field regressed from time series data generated by a high-dimensional vector field 
$F \colon \mathbb{R}^N \to \mathbb{R}^N$, with $n \ll N$. (Indeed, in the applications 
discussed below, $F$ comes from an infinite-dimensional problem. See Section \ref{sec:regression}
for a more formal description of this regression). One natural diagnostic is to measure how 
long trajectories of $f$ track corresponding observed trajectories of $F$, for initial conditions 
not used in the regression. However, sensitivity to initial conditions and 
modeling errors severely limits the time horizon over which we expect individual trajectories 
to remain correlated. More importantly, accurate finite-time tracking does not, by itself, 
establish whether recurrent invariant structures of the surrogate $f$ correspond to genuine 
invariant structures of the original system $F$. So, can recurrent invariant structures identified 
in the surrogate be used to establish the existence of corresponding invariant structures 
in the original, potentially infinite-dimensional, dynamical system?

The general question is far beyond the scope of any single investigation. Here we propose 
a diagnostic for one particular kind of recurrence: periodic orbits. More precisely, we use 
harmonic balance methods to numerically compute Fourier expansions of periodic orbits in the 
low-dimensional surrogate, and with these obtain approximate solutions in the original system. 
We then apply a-posteriori computer-assisted proofs to certify the existence of a true periodic 
solution nearby. When successful, this procedure provides rigorous bounds on the distance 
between the surrogate-generated approximation and a genuine invariant solution of the original, 
high-dimensional system. We show that such certifications can be obtained without the use of 
any iterative correction scheme in the high-dimensional system. Thus, the existence proof itself 
becomes a quantitative diagnostic of how well the surrogate reproduces an invariant structure 
in the original system.

For surrogate construction, we use a technique we refer to as dense polynomial surrogate 
reconstruction (DPSR). This consists of fitting continuous-time polynomial vector fields 
by least squares or ridge regression, retaining a complete polynomial library in prescribed 
observation coordinates. Such polynomial reconstruction methods have a substantial 
history \cite{gouesbet1994}, and the same general regression framework underlies more recent 
sparse identification methods \cite{brunton2016}. An important distinction is that DPSR 
reconstructs a continuous-time vector field, rather than a discrete-time evolution map, 
as is common in classical nonlinear time-series prediction \cite{farmer1987,casdagli1989} 
and modern reservoir-computing approaches \cite{pathak2017,gauthier2021}. This choice is 
motivated by our interest in computing invariant structures using established techniques 
from computational dynamics. The resulting ordinary differential equations are readily 
differentiated and studied using standard methods for locating equilibria, periodic orbits,
and other invariant sets. We make no claim of novelty for the regression procedure, which 
is chosen for its simplicity. The certification strategy itself is not tied to this or any
particular regression method: other techniques could be used, provided that suitable 
approximate solutions of the original system can be reconstructed from the learned dynamics.

We illustrate the procedure by taking simulated time series data from the Mackey--Glass 
delay differential equation and the Kuramoto--Sivashinsky (KS) partial differential equation, 
both of which generate infinite-dimensional dynamical systems, and both of which
are well known to exhibit complex dynamics.  For Mackey--Glass, we regress a vector field using 
delay coordinates, motivated by classical embedding theory \cite{takens1981,sauer1991}, and 
finite-difference estimates of derivatives. For KS, we use spatial coordinates and 
derivative data obtained from spectral or Galerkin approximations of the governing equations. 
Held-out forecasts, long-time integrations, equilibrium spectra, finite-time Lyapunov estimates, 
and attractor projections provide complementary numerical assessments. Such dynamical diagnostics 
have also been used to evaluate machine-learned models of chaotic systems \cite{pathak2017}. 
Related work has used machine-learned representations to identify invariant solutions 
organizing turbulent dynamics \cite{page2021}.

The rigorous component of our approach uses established methods for computer-assisted proofs 
of invariant objects. These techniques combine bounds on the residual, an approximate inverse, 
and nonlinear terms to establish existence through contraction arguments in suitable function 
spaces \cite{framework2017}. Such methods have previously been applied to periodic solutions of 
KS \cite{figueras2017} and delay equations, including Mackey--Glass \cite{berg2022}. Our 
contribution is to use these methods not simply to establish existence, but to assess the ability 
of a learned surrogate to produce dynamically meaningful invariant structures. To this end, 
we freeze the surrogate-generated approximation before the rigorous validation step, ensuring 
that the resulting existence theorem provides an error bound directly for the surrogate's output.

The Mackey--Glass and KS examples demonstrate that this strategy can successfully certify 
periodic solutions of infinite-dimensional systems using candidates constructed from 
low-dimensional surrogates, without Newton corrections in the original equations. 
Subsequent refinement using the original equations produces more 
accurate approximations of the certified periodic solutions, but is not needed for the 
initial existence results.

These results provide rigorous, local information about the correspondence between recurrent 
structures of a surrogate and those of the original system. They do not establish global 
equivalence of the dynamics, orbital stability, or agreement of long-time statistical 
properties. Conversely, failure of an a-posteriori existence test does not necessarily indicate 
that a surrogate is inaccurate, since the available proof bounds may be insufficient. 
Nevertheless, successful certification provides a mathematically precise assessment of 
structural fidelity that cannot be obtained from finite-time trajectory comparisons alone.

Our first example is the Mackey--Glass equation \cite{mackey1977}, a scalar
delay differential equation originally proposed as a model of physiological
control. Because its state at each time is a history segment of length equal
to the delay, its phase space is infinite-dimensional, and at the standard
parameter values its dynamics are chaotic. We observe only a single scalar time
series, embed it in six delay coordinates, and fit a polynomial vector field in
these coordinates. A periodic candidate obtained from the fitted model is then
validated directly in the delay equation with all physical parameters fixed,
using Fourier methods for delay equations in the spirit of \citet{berg2022}.
The resulting certificate encloses a true periodic solution within $6\times10^{-4}$
of the surrogate export, with no Newton corrections in the delay equation.

Our second example is the Kuramoto--Sivashinsky equation
$u_t+uu_x+u_{xx}+u_{xxxx}=0$ on an interval of length $22$, a standard model of
spatiotemporal chaos. Here the long-wave instability $u_{xx}$ is balanced by the
strong damping $u_{xxxx}$ of short waves, so that the dynamics are effectively
low-dimensional despite the infinite-dimensional phase space
\cite{foias1988}. We use two settings. With periodic boundary conditions, a
quadratic model in the values at 14 spatial nodes is assessed by forecasts,
attractor projections, and equilibria. With Neumann conditions for the potential,
the velocity is a sine series, and the computer-assisted framework of
\citet{figueras2017} applies. A quadratic surrogate in the sine coefficients yields
a periodic candidate that is certified within $10^{-5}$ of the export, again without
Newton corrections. Since that surrogate is trained on exact Galerkin derivatives,
we also learn Neumann KS surrogates from time series recorded at fixed spatial
sites, as in a network of weather stations. These station surrogates produce
candidates whose error decreases with the number of sensors and which converge
to true periodic orbits after a few Newton corrections. Beyond about 24 sensors,
however, the regression no longer identifies the damping of the finest resolved
modes, and the fitted models become unstable. The certification-based diagnostic
thus records both a success and a limitation of the polynomial surrogates.

The remainder of the paper is organized as follows. Section~\ref{sec:regression} describes 
the regression construction and numerical assessment protocols. Section~\ref{sec:results} 
presents forecasting results and numerical dynamical diagnostics for the learned surrogate 
models. Section~\ref{sec:assessment} develops the proposed candidate-based validation 
procedure, while Sections~\ref{sec:cap} and~\ref{sec:cap-results} describe the 
computer-assisted proof framework and certification results. Conclusions and possible 
extensions are discussed in Section~\ref{sec:discussion}.

\section{Polynomial regression and surrogate construction}
\label{sec:regression}

We now describe the regression procedure used to construct the low-dimensional 
surrogate models. The basic idea is to fit an explicit polynomial vector field to 
time series data coming from observed states of a higher dimensional system. 
While the regression algorithm is common to all our 
experiments, the observation coordinates and derivative targets differ between 
applications.  We consider three common situations.  The first is that we have 
only a single scalar time series, and a higher dimensional phase is reconstructed
through the use of delay embedding.  In the second case we have access to 
time series data coming from a spatially extended system, 
measured at a number of different ``sites.''  Finally, we consider the case that
the spatially extended system is described using a Galerkin projection, and 
suppose that we have access to data from a finite number of modes.

\subsection{Observation coordinates and polynomial library for a scalar time series}
\label{sec:delay-coordinates}

Given scalar time series data from 
an observed state variable of a higher dimensional system,  
we embed the time series into a higher dimensional phase space
using a classic delay embedding technique.  That is, given samples $y_j=x(t_0+jh)$, 
we choose an integer lag $\ell\geq 1$, corresponding delay $\Delta=\ell h$, 
and embedding dimension $d$. For indices at which all delayed observations 
are available, define
\begin{equation}
X_j=(y_j,y_{j-\ell},\ldots,y_{j-(d-1)\ell})^\top\in\R^d.
\label{eq:sampled-delay-state}
\end{equation}
This construction is motivated by classical embedding 
theory \cite{takens1981,sauer1991}, although the regression procedure itself 
requires only a collection of state--derivative pairs.

After re-indexing the $M$ available states, the data matrix $\mathbf X\in\R^{M\times d}$ 
has rows $X_j^\top$. We estimate derivatives using second-order finite differences:
\begin{equation}
\begin{aligned}
\dot X_j&=\frac{X_{j+1}-X_{j-1}}{2h},
&&1\leq j\leq M-2,\\
\dot X_0&=\frac{-3X_0+4X_1-X_2}{2h},\\
\dot X_{M-1}&=\frac{3X_{M-1}-4X_{M-2}+X_{M-3}}{2h}.
\end{aligned}
\label{eq:derivatives}
\end{equation}
The resulting derivative estimates form a matrix $\Xdot\in\R^{M\times d}$.

We seek a polynomial vector field $f_C:\R^d\to\R^d$ approximating the observed derivatives. 
We choose a maximum degree $p$, and introduce the complete polynomial library
\begin{equation}
\mathcal A_{d,p}
=\{a\in\Nzero^d:|a|\leq p\},
\qquad
\theta_a(X)=\prod_{k=1}^{d}X_k^{a_k},
\qquad
q=\binom{d+p}{p},
\label{eq:polynomial-library}
\end{equation}
where $|a|=\sum_{k=1}^d a_k$. Ordering the multi-indices with $a^{(0)}=0$, 
define the feature vector
\[
\Theta_p(X)
=\bigl(\theta_{a^{(r)}}(X)\bigr)_{r=0}^{q-1}
\]
and the design matrix
\[
\mathbf\Theta_{j,r}
=\theta_{a^{(r)}}(X_j).
\]
Our surrogate is the autonomous polynomial differential equation
\[
\dot X=f_C(X)=C\Theta_p(X),
\qquad C\in\R^{d\times q},
\]
whose coefficients are determined by regression.

This is the standard polynomial-library formulation used in data-driven equation 
discovery \cite{brunton2016}. In contrast to sparse identification methods, we initially 
retain every monomial through degree $p$. We also emphasize that an autonomous vector 
field fitted in delay coordinates need not preserve exact delayed-copy relations 
among its components, an issue that becomes relevant when constructing 
periodic-orbit candidates.

\subsection{Normalized ridge regression and model selection}
\label{sec:ridge-regression}

The monomials in the library may have dramatically different magnitudes on the training 
data. To improve numerical conditioning and control the influence of higher-degree 
terms, we normalize the features and use degree-weighted ridge regularization.
That is, for each feature, define its root-mean-square magnitude
\begin{equation}
s_r=\left(M^{-1}\sum_j\mathbf\Theta_{j,r}^2\right)^{1/2},
\qquad
\Thetahat_{j,r}=\mathbf\Theta_{j,r}/s_r,
\label{eq:feature-normalization}
\end{equation}
with $s_r$ replaced by one whenever $s_r<10^{-14}$. The regularization weights are
\[
w_0=0,\qquad
w_r=\max\{|a^{(r)}|,1\}^{\kappa},
\quad r\geq1,
\]
where $\kappa\geq0$. In particular, the constant term is not penalized.

For a regularization parameter $\lambda\geq0$, we determine the normalized 
coefficients from
\begin{equation}
\widehat C\in
\mathop{\rm argmin}_{B\in\R^{d\times q}}
\left\{
\|\Thetahat B^\top-\Xdot\|_F^2
+
\lambda\sum_r w_r\|B_{:,r}\|_2^2
\right\},
\qquad
C_{:,r}=\widehat C_{:,r}/s_r.
\label{eq:ridge-objective}
\end{equation}
This is a degree-weighted version of classical ridge regression \cite{hoerl1970}.
We solve the regression problem using rank-revealing least squares applied 
to the augmented design matrix
\[
\begin{bmatrix}
\Thetahat\\
\sqrt{\lambda}\diag(\sqrt w)
\end{bmatrix},
\]
with corresponding target
\[
\begin{bmatrix}
\Xdot\\
0
\end{bmatrix}.
\]
This avoids forming normal equations, which can amplify conditioning difficulties. 
When $\lambda=0$, a rank-deficient library is handled using the minimum-norm 
solution in the chosen feature scaling.

A small derivative residual does not necessarily imply that the fitted system has 
accurate long-time dynamics. We therefore select the parameters $(\Delta,d,p,\lambda)$ 
using validation data disjoint from the training data, considering forecast errors, 
numerical boundedness, and qualitative diagnostics such as amplitude, variance, 
autocorrelation, and spectral content. Separate test data are reserved for 
subsequent assessment.

\begin{remark}[Dense polynomial regression and relation to SINDy]
\label{rem:sindy}

Our regression approach is closely related to sparse identification of nonlinear 
dynamics (SINDy) \cite{brunton2016}. Both methods approximate a vector field using a 
prescribed library of candidate functions. The main distinction is that SINDy seeks a 
small active set of terms, often using sequentially thresholded least squares, 
whereas our baseline models retain the complete polynomial library.
Complete-library polynomial regression predates SINDy \cite{crutchfield1987,gouesbet1994,brown1994}. Related work treats unobserved variables \cite{breeden1990}, fits differential equations to chaotic data by boundary-value methods \cite{baake1992}, identifies spatiotemporal systems and partial differential equations from data \cite{voss1998,bar1999,rudy2017}, and searches over symbolic expressions or sparse expansions \cite{bongard2007,schmidt2009,wang2011}; see \cite{aguirre2009} for a review. In the present work, ridge 
regularization controls coefficient magnitudes without imposing a sparse support. 
Pruning and library augmentation are considered separately as numerical ablations, 
rather than as part of the baseline fitting procedure.
\end{remark}

\subsection{Spatial coordinates and KS surrogates}

When we consider PDEs, we can take observation data at spatial sites 
rather than delay embedding reconstructions. We consider two settings, with different 
boundary conditions and distinct roles in the experiments.
For periodic KS on $[0,L)$, we represent the state by values at $d$ equally spaced 
spatial nodes:
\[
X(t)=(u(x_1,t),\ldots,u(x_d,t))^\top,
\qquad x_j=\frac{(j-1)L}{d}.
\]
The derivative targets are computed from the reference spectral discretization, 
using dealiased spectral differentiation at the training snapshots.

The main forecasting surrogate retains all quadratic monomials and is obtained by 
unscaled, unregularized dense least squares, without coefficient thresholding. 
If the design matrix is rank deficient, the SVD-based solver selects its minimum-norm 
coefficient vector in these unscaled coordinates. Thus, unlike the delayed scalar
case, this fit uses neither feature normalization nor ridge 
regularization. The numerical rank and singular values are recorded to characterize 
the regression problem.

For KS with Neumann potential boundary conditions, we instead work with a 
truncated sine representation of the velocity $u=h_x$:
\begin{equation}
u(x,t)=\sum_{k=1}^{K}a_k(t)\sin(k\pi x/\ell),
\qquad
\dot a=c+La+Q(a,a).
\label{eq:sine-coordinates}
\end{equation}
Here $a\in\R^K$ is the vector of sine coefficients, and the right-hand side is the 
corresponding Galerkin system, consisting of constant, linear, and quadratic 
terms. The coefficients $a_k$ and the nodal values at
\[
x_j=\frac{j\ell}{K+1},
\qquad 1\leq j\leq K,
\]
are related by an invertible discrete sine transform.

Unlike the previous two experiments, the Neumann KS training data are not generated 
from a time series. Instead, we sample independent states in coefficient space and
evaluate their derivatives using the known Galerkin equations. Specifically, we take
\[
a_k=\sigma_k z_k,
\qquad
\sigma_k=\left(1+\frac{k-1}{4}\right)^{-2},
\]
where the $z_k$ are independent standard normal random variables.

The quadratic structure of the Galerkin vector field allows us to separate its odd
and even parts by evaluating the derivative target $g$ at paired states $a$ and $-a$:
\begin{equation}
\frac12\bigl(g(a)-g(-a)\bigr)=La,
\qquad
\frac12\bigl(g(a)+g(-a)\bigr)=c+Q(a,a).
\end{equation}
We fit these two parts separately by rank-revealing least squares in scaled 
coordinates. With
\[
q_e=1+\frac{K(K+1)}{2}
\]
even polynomial features, we generate $2q_e$ independent states and their negatives, 
giving $4q_e$ training snapshots. An additional 512 independently sampled states are 
used to assess the derivative fit.

We also enforce the reflection symmetry of the Neumann problem. In coefficient 
space, define
\[
(Ra)_k=(-1)^k a_k.
\]
Reflection equivariance is imposed by averaging the fitted vector field with its 
reflection conjugate. The resulting polynomial model respects this symmetry by 
construction. No selected periodic orbit is included in the training ensemble.

An important feature of this experiment is that the Galerkin right-hand side is itself 
quadratic. Consequently, the target vector field belongs to the chosen polynomial 
library. In exact arithmetic, if the design matrix has full column rank, regression 
recovers this vector field exactly. Indeed, if
\[
D=\Theta C_*^\top,
\qquad \ker\Theta=\{0\},
\]
then $C_*$ is the unique zero-residual coefficient matrix.

This observation explains the near-machine-precision derivative errors obtained in the 
Neumann KS regression. Such accuracy demonstrates recovery of the finite-dimensional 
Galerkin operator, but does not by itself establish accuracy with respect to the full 
PDE. The latter question is addressed by the rigorous certification procedure.

\paragraph{Regression implementation.}
The shared polynomial library, regression solvers, and vector-field evaluation 
routines are provided in \path{dpsr/polynomial.py}, with experiment configurations 
in \path{notebooks/} in the \href{https://github.com/Aloh5/learned-dynamics-cap}{project's GitHub repository}. The Neumann sampling and symmetry routines are implemented 
in \path{CAP/Kuramoto_PO/surrogate_neumann.py}.

\subsection{Station observations for Neumann KS}
\label{sec:stations}
The Galerkin-supervised Neumann fit above is a controlled best case. To test
the method on data of the kind available in practice, we also learn a
Neumann KS surrogate from synchronized time series at $M$ fixed spatial sites
(``weather stations''). In this experiment the Galerkin field is used only to
generate the data and, after a candidate has been frozen, to evaluate it.

\paragraph{Sites.}
We place the stations on the regular grid
\begin{equation}
 x_j=\frac{j\ell}{M+1},\qquad 1\leq j\leq M,
 \label{eq:station-sites}
\end{equation}
which excludes the endpoints, where $u=0$ for every state. On this grid the
nodal values and the first $M$ sine coefficients are related by the discrete
sine transform (DST)
\begin{equation}
 Y_j=\sum_{k=1}^{M}a_k\sin\frac{jk\pi}{M+1},
 \qquad
 a_k=\frac{2}{M+1}\sum_{j=1}^{M}Y_j\sin\frac{jk\pi}{M+1},
 \label{eq:dst}
\end{equation}
since the matrix $S_{jk}=\sin(jk\pi/(M+1))$ satisfies $S^\top S=\tfrac{M+1}{2}I$.
This choice makes reconstruction exact for profiles in the span of the first $M$
sine modes. Modes $k>M$ alias onto retained modes, so the reconstruction error
decreases with the spectral decay of the solution. We also tested irregular
sites, drawn uniformly at random in $(0,\ell)$ with minimum spacing
$\ell/(4(M+1))$, combined with reconstruction by a natural cubic spline through
$(0,0),(x_j,Y_j),(\ell,0)$ followed by $L^2$ projection onto sine modes. With
exact orbit data, this route was 1--3 orders of magnitude less accurate than
\eqref{eq:dst}, and every irregular-site surrogate with $M\leq24$ diverged.
We therefore report the regular grid.

\paragraph{Data.}
The reference trajectory is the 48-mode sine--Galerkin system integrated by
ETDRK4 with step $0.025$ from a fixed seeded initial state. After a burn-in of
300 time units, we record $Y(t)=(u(x_1,t),\ldots,u(x_M,t))^\top$ every $h=0.1$
for 3,000 time units, giving 30,001 snapshots. All sensor counts use the same
trajectory. No noise is added.

\paragraph{Regression.}
Time derivatives are estimated from the samples either by the second-order
differences \eqref{eq:derivatives} (FD2) or by fourth-order central
differences (FD4), with second-order one-sided formulas at the two ends. We fit
a complete quadratic $\dot Y=f(Y)$ in the nodal coordinates by the normalized,
degree-weighted ridge regression of Section~\ref{sec:ridge-regression}, using
the first 70\% of the snapshots. For $\lambda\in\{0,10^{-10},10^{-8},10^{-6}\}$
we compute the relative derivative error on the remaining 30\%. Following the
model-selection rule of that section, we take the $\lambda$ with the smallest
held-out error whose RK4 rollout (step $0.005$, or $0.002$ for $M\geq36$) stays
finite and bounded by $10^3$ for 2,300 time units, starting from the first
held-out snapshot. If every $\lambda$ fails, the sensor count is recorded as
divergent. No symmetry is imposed on the nodal fit.

\paragraph{Periodic orbits of the surrogate.}
After discarding the first 300 time units of the rollout, we scan for close
returns with periods in $[20,120]$ and use the eight best as seeds. Each
segment is resampled to 33 equally spaced phases. We then solve the
harmonic-balance (Fourier collocation) equations
\begin{equation}
 \omega\,D\,\mathbf Y-f(\mathbf Y)=0,\qquad
 \sum_i\langle \mathbf Y_i-\mathbf Y^{\rm seed}_i,\,(D\mathbf Y^{\rm seed})_i\rangle=0,
 \label{eq:station-hb}
\end{equation}
for the nodal values $\mathbf Y$ at the phase nodes and the frequency $\omega$,
where $D$ is the Fourier differentiation matrix in phase. We use 16 and then 32
temporal harmonics, with a Gauss--Newton solver. We accept an orbit if its
collocation RMS residual is below $10^{-9}$, its amplitude exceeds $10^{-2}$,
and its period differs from those already found.

\paragraph{Candidate and evaluation.}
The temporal Fourier coefficients $b_n\in\mathbb{C}^M$ of an accepted orbit are
mapped to spatial sine coefficients by \eqref{eq:dst}, applied to each $n$. This
gives a space--time candidate with 32 temporal harmonics and $M$ sine modes,
which is then frozen. Only after freezing do we evaluate it in the true
Galerkin collocation equation, with 32 temporal harmonics and 48 sine modes.
The \emph{initial defect} is the maximum collocation residual of the frozen
candidate. We then apply plain Newton iterations to the phase-augmented true
equation, starting from the candidate, until the maximum residual is below
$10^{-10}$. We record the number of corrections and the distance in the
norm~\eqref{eq:ksnorm} between the frozen candidate and the Newton limit.
These are floating-point diagnostics. Certification is discussed in
Section~\ref{sec:station-results}.

\subsection{Assessment protocols}
\label{sec:protocols}

We assess the fitted surrogates using both finite-time forecasting and numerical diagnostics of their longer-time dynamics.

For forecasting, let $e_j$ denote the scalar prediction error for Mackey--Glass or the spatial root-mean-square error for KS. We define the forecast horizon by the first crossing of a cumulative squared-error threshold:
\begin{equation}
H_\varepsilon
=\min\left\{
jh:h\sum_{i=0}^j e_i^2>\varepsilon
\right\},
\qquad
\varepsilon_{\rm MG}=0.05,
\qquad
\varepsilon_{\rm KS}=0.5.
\label{eq:horizon}
\end{equation}
If the threshold is not crossed, the horizon is set equal to the forecast duration.

We also report normalized root-mean-square error (NRMSE), obtained by dividing the full-window RMSE by the standard deviation of the reference rollout. Unlike the forecast horizon, this statistic includes errors accumulated after any threshold crossing. Since the two systems have different state scales and error thresholds, their forecast horizons are interpreted separately rather than compared directly.

Beyond forecasting, we examine the dynamics generated by the learned vector fields. Numerical diagnostics include long-time boundedness, empirical distributions, variance, autocorrelation, and spectral content. We also estimate trajectory separation rates using repeated perturbation renormalization, providing finite-time approximations of Lyapunov exponents. These tests complement short-time forecasting but do not establish rigorous properties of invariant sets.

\paragraph{Reference data.}
For the main Mackey--Glass forecasting experiment, the reference trajectory is integrated using Euler's method with time step $0.1$, while the surrogate is integrated using RK4 with the same step. We accept surrogate trajectories whose coordinates remain finite and within $[0.15,1.65]$.

The training, validation, test, and reference-statistics windows are
\[
[117,417),\quad [700,1300),\quad
[1500,2500),\quad [3500,4500),
\]
respectively. The reference-statistics window is used only for assessment and is not included in model selection.

Forecasts starting at different times within a single trajectory may have overlapping evaluation windows. Consequently, the reported medians and ranges are descriptive statistics of these experiments, not confidence estimates based on independent realizations.

\section{Review of forecasting and numerical assessment results}
\label{sec:results}
This section is not a methodological contribution. Fitting polynomial or
other basis-function models in delay or spatial coordinates and assessing
them by held-out forecasts, attractor statistics, and Lyapunov estimates
are established practice
\cite{crutchfield1987,farmer1987,casdagli1989,brown1994,abarbanel1993,pathak2017,pathak2018,gauthier2021};
Mackey--Glass and KS are standard benchmarks in that literature. We report
these results to document the surrogates that supply the candidates of
Sections~\ref{sec:assessment}--\ref{sec:cap-results} and to place the later
certificates alongside conventional evidence.
Three settings support distinct claims: Mackey--Glass assesses delayed-scalar
forecasting and local dynamics; periodic KS assesses a 14-node spectral
semidiscretization; Neumann-potential KS supplies an orbit candidate for rigorous
validation. The two KS problems differ in boundary conditions and coordinates.
Forecasts use the assessment protocols in Section~\ref{sec:protocols}.

\subsection{Mackey--Glass: prediction and local dynamics}
\label{sec:mg-results}
\begin{figure}[htb]
\centering
\includegraphics[width=\textwidth]{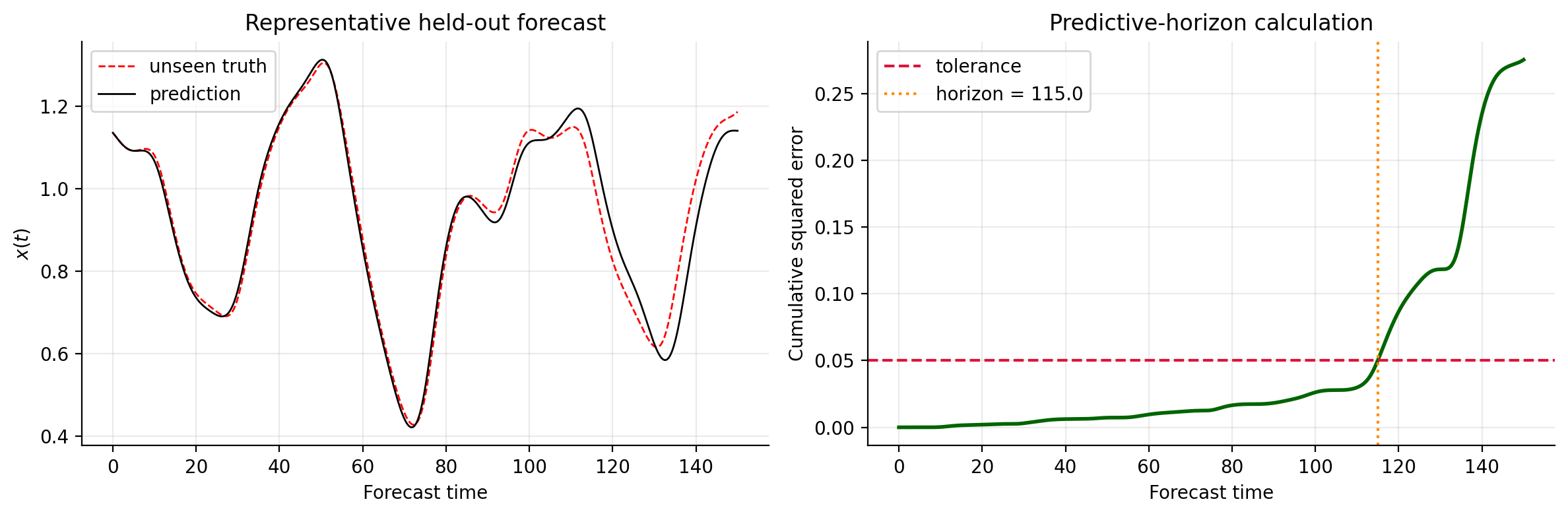}
\caption{Main cubic Mackey--Glass forecast on a disjoint held-out window.}
\label{fig:mgforecast}
\end{figure}

The reference equation, with constant history $x(t)=1.2$, is \cite{mackey1977}
\begin{equation}
 x'(t)=\frac{0.2x(t-17)}{1+x(t-17)^{10}}-0.1x(t).
 \label{eq:mgphysical}
\end{equation}
All 12 forecasts of the main cubic model complete 150 time units within the
numerical bounds (Figure~\ref{fig:mgforecast}). Horizon has mean $99.34$, median
$99.7$, and range $67.1$--$140.2$; full-window NRMSE has mean $0.1676$, median
$0.1715$, and range $0.0862$--$0.2994$. These starts sample overlapping windows
of one held-out segment.

\paragraph{Numerical evidence of chaos.}
Three integrations reach 20,000 time units while remaining within the prescribed
numerical bounds. For each initial condition, a two-trajectory Benettin-style
calculation discards a 5,000-unit transient, accumulates logarithmic separation
growth for a further 2,000 units, and renormalizes a perturbation of size
$10^{-7}$ every unit of time, using RK4 at step $0.1$. All three resulting
maximal Lyapunov estimates are positive (Table~\ref{tab:chaos}). These are
finite-time numerical estimates, not validated asymptotic exponents.

\begin{table}[htb]
\centering\small
\caption{Long-run diagnostics of the main cubic Mackey--Glass surrogate.
Statistics and peak counts use the final 5,000 units of each 20,000-unit run;
Lyapunov accumulation uses the separate post-transient protocol in the text.
Standard deviations are relative to the reference-statistics window.
Distinct peak heights are counted after rounding to $10^{-4}$.}
\label{tab:chaos}
\begin{tabular}{crrrr}
\toprule
Run & Lyapunov estimate & Std.\ ratio & Peaks & Distinct heights\\
\midrule
1 & 0.006405 & 1.0103 & 158 & 153\\
2 & 0.007114 & 1.0136 & 163 & 156\\
3 & 0.005867 & 1.0133 & 151 & 144\\
\bottomrule
\end{tabular}
\end{table}

\begin{figure}[htb]
\centering
\includegraphics[width=\textwidth]{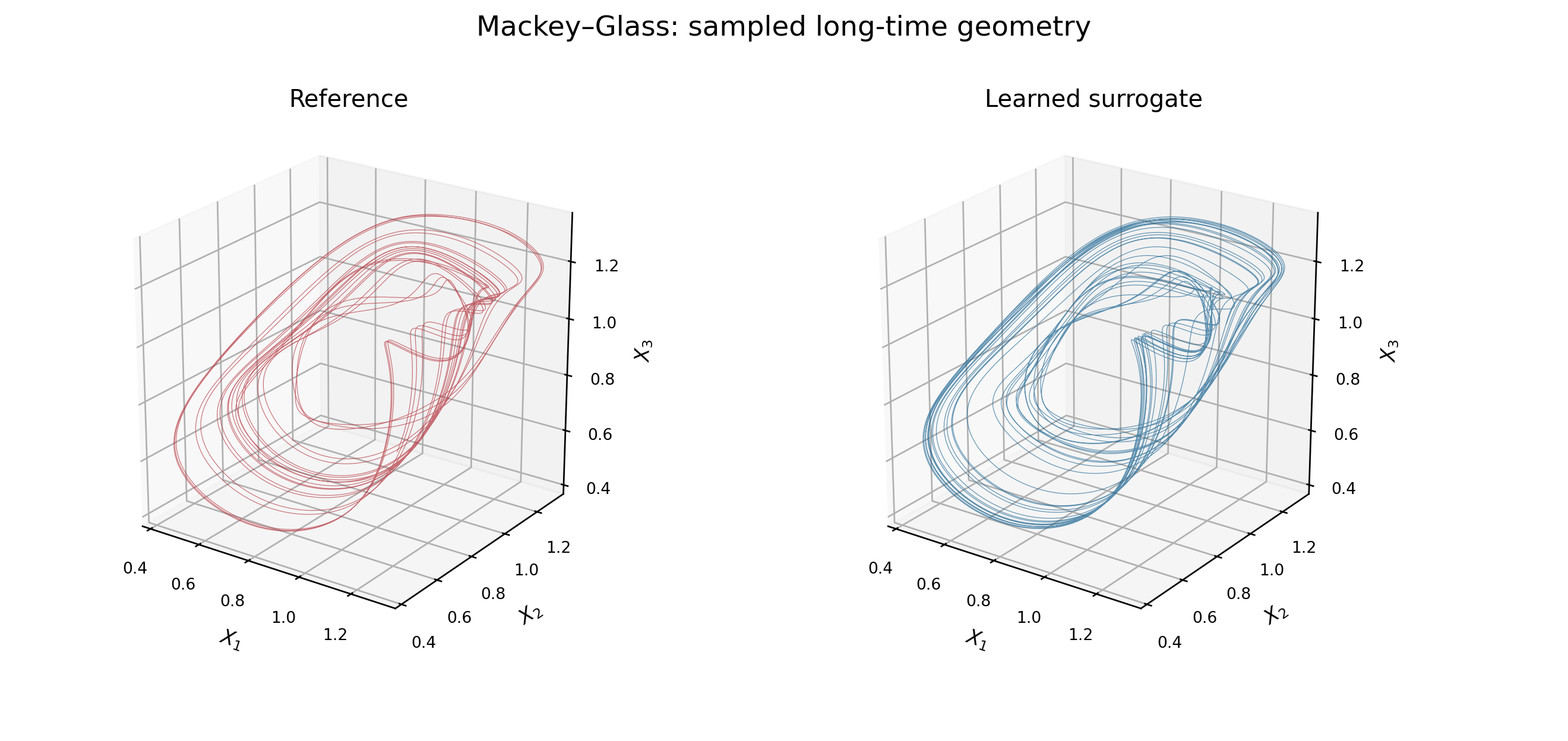}
\caption{Mackey--Glass attractor comparison in the first three observation
coordinates, with identical axis limits and viewing angles. Left: reference
Euler data at step $0.1$ on $[3500,4500)$, with
$(X_1,X_2,X_3)=(x(t),x(t-4.5),x(t-9))$. Right: the stored main cubic surrogate,
integrated with RK4 at step $0.1$ for 2,000 units from the reference state at
$t=1500$, with the first 500 units discarded. The surrogate coordinates
need not remain exact delayed copies. These finite trajectory projections
compare sampled geometry; they do not prove equality of attractors.}
\label{fig:mg-attractor-comparison}
\end{figure}

\begin{figure}[htb]
\centering
\includegraphics[width=\textwidth]{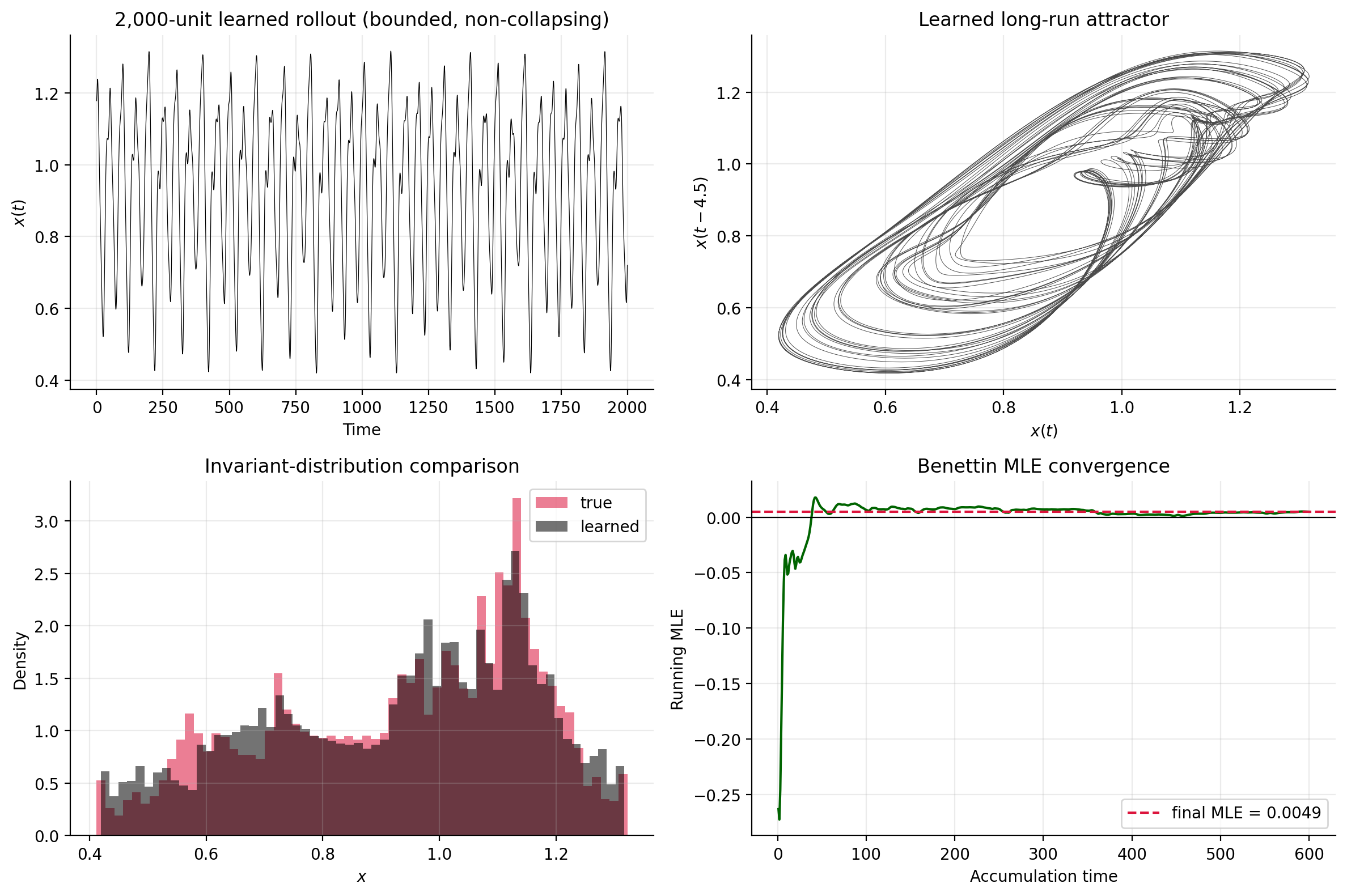}
\caption{Qualitative Mackey--Glass diagnostics from the saved cubic-model
experiment: a 2,000-unit rollout, its first-two-coordinate projection, scalar
histograms, and a running Lyapunov estimate. The histogram discards the first
500 units of the learned rollout. The displayed Lyapunov value, approximately
$0.0049$, uses a separate 200-unit transient and 600-unit accumulation and is
not one of the longer estimates in Table~\ref{tab:chaos}. ``Attractor'' and
``invariant-distribution'' in the panel titles refer to numerical projections
and empirical histograms, not rigorous invariant-set or measure claims.}
\label{fig:mgattractor}
\end{figure}

\paragraph{Qualitative attractor diagnostics.}
Figure~\ref{fig:mg-attractor-comparison} compares long-time reference and
surrogate geometry directly in the first three observation coordinates.
Figure~\ref{fig:mgattractor} shows the existing 2,000-unit learned rollout,
a projection onto the first two surrogate coordinates, and empirical scalar
distributions for the learned and reference trajectories. The rollout remains
irregular and explores a folded band of states rather than visibly collapsing
onto one simple loop. The distributions have similar range and broad shape,
while their local differences remain visible. The plotted coordinate labels
use the delay-coordinate convention; a free surrogate trajectory need not
satisfy exact delayed-copy compatibility.

The separate 20,000-unit runs provide a longer check: their standard deviations
are within approximately $1.4\%$ of the reference value, and most late peak
heights remain distinct (Table~\ref{tab:chaos}). Positive finite-time separation
growth, persistent amplitude variation, and the projected geometry together
support chaotic-like dynamics of the cubic surrogate over the sampled intervals.
These are qualitative numerical checks of attractor behavior. They do not
prove a chaotic invariant set, equality of invariant measures, a conjugacy, or
global boundedness. The periodic-orbit theorems below establish different,
local existence statements.

A search from 279 starts finds 11 equilibria with field residual below $10^{-8}$.
The root nearest $\mathbf1\in\R^6$ has error $2.59\times10^{-2}$ and unstable
eigenvalues $0.03962\pm0.13938i$, versus the DDE's leading roots
$0.04129\pm0.13907i$ at $x=1$; both have two unstable directions. This root is
also closest to the sampled training trajectory (distance $0.2167$, versus
$0.3246$ for the next closest). The other ten roots have inconsistent delay
coordinates and lie at least $0.5735$ from the constant-equilibrium embeddings
$-\mathbf1,\mathbf0,\mathbf1$. Thus local instability is approximated despite
additional surrogate equilibria. The search is finite, and proximity to training
data does not establish attractor membership.

The degree-five orbit-refinement model has a narrower role: global-rollout
tests cross the prescribed escape threshold near $57.91$ from a reference-attractor
state and $2.81$ from a perturbed state. Local candidate refinement therefore
does not imply improved global forecasting.

\subsection{Periodic KS: prediction and invariant structures}
\label{sec:periodic-ks-results}
The periodic reference is
\begin{equation}
 u_t=-uu_x-u_{xx}-u_{xxxx},\qquad u(x+L,t)=u(x,t),\quad L=22.
 \label{eq:ksperiodic}
\end{equation}
The main quadratic model uses 14 nodes, 120 features, and 1,680 coefficients.
ETDRK4 at step $0.025$, sampled every $0.1$, generates a 100-unit burn-in,
70-unit training window (701 snapshots), and 100-unit forecast (1,001 snapshots)
\cite{kassam2005}. The shared endpoint initializes the forecast. The training
library has numerical rank 106, and the relative training derivative residual
is $5.64\times10^{-14}$.

The fixed quadratic model was evaluated without refitting on 100 independently
drawn smooth initial fields with the training spatial mean, each followed
by a 100-unit burn-in and a 100-unit forecast. All 100 rollouts completed;
99 remained below the cumulative spatial-error threshold
$\varepsilon_{\rm KS}=0.5$ throughout, while one crossed at $74.2$.
Full-window NRMSE had median $4.97\times10^{-4}$ and mean
$7.08\times10^{-3}$, including the largest error of $0.5533$.
The 99 capped horizons establish $H_{0.5}\geq100$.
Reference trajectories used ETDRK4 step $0.0125$ and surrogate forecasts used
RK4 step $0.1$. These results concern the fixed-mean, 14-node
semidiscretization; fully resolved PDE forecasting requires a
spatial-convergence study.

\begin{figure}[htb]
\centering
\includegraphics[width=\textwidth]{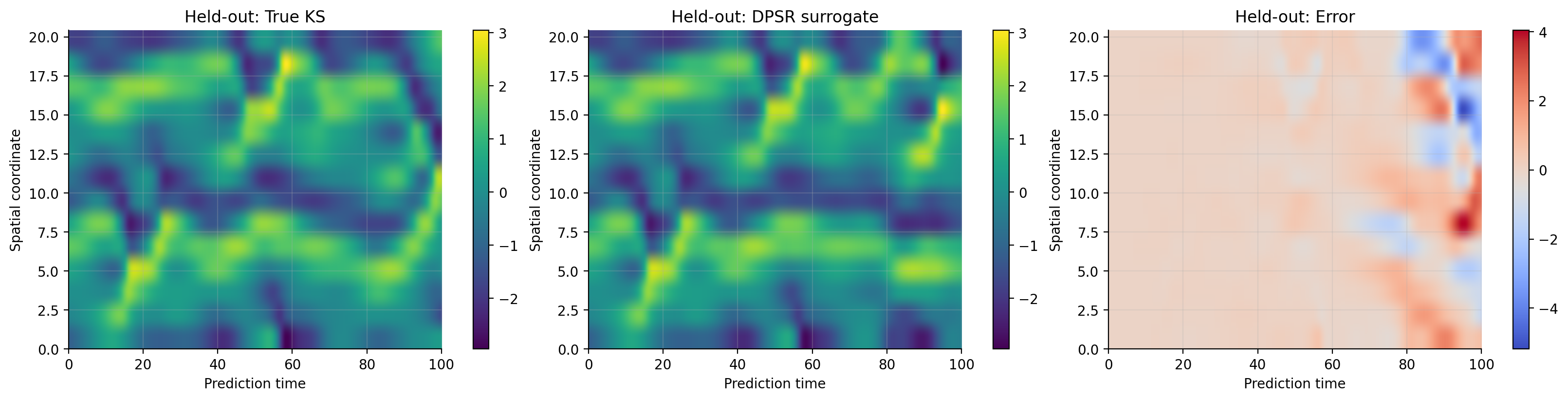}
\caption{Periodic KS: reference, surrogate, and error over the 100-unit forecast.}
\label{fig:ksforecast}
\end{figure}

\begin{figure}[htb]
\centering
\includegraphics[width=\textwidth]{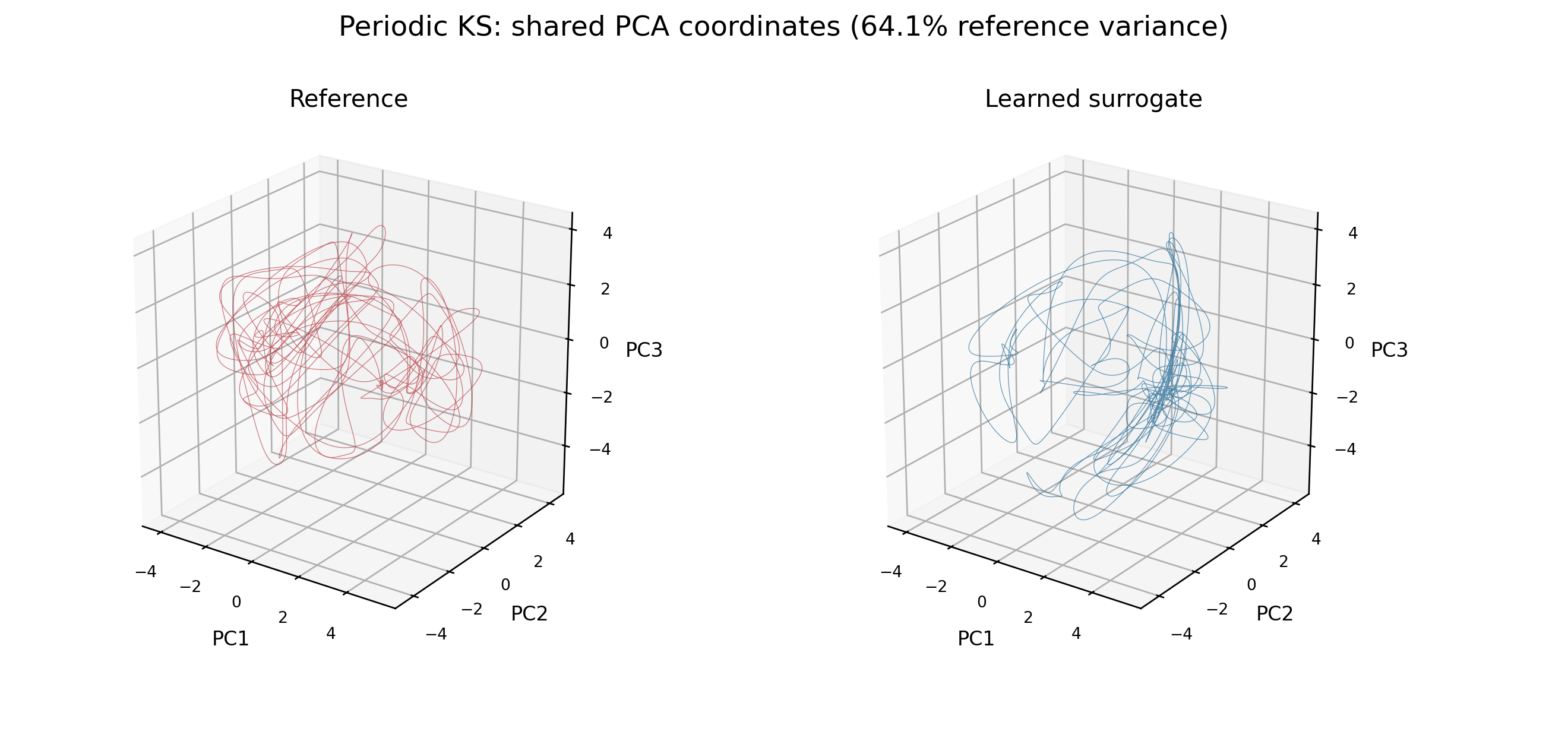}
\caption{Periodic KS attractor comparison for $L=22$ and 14 nodes.
Reference: ETDRK4 with step $0.025$, sampled every $0.1$.
Learned surrogate: the stored quadratic forecasting model, integrated with
DOP853, relative tolerance $10^{-9}$, absolute tolerance $10^{-11}$, maximum
step $0.1$, and output spacing $0.1$. Both panels show the retained interval
$[200,1000]$ in the same three-dimensional PCA projection; the title reports
the reference variance captured. No model is refitted for this figure.
The comparison concerns a spectral semidiscretization, not a
spatially converged PDE attractor.}
\label{fig:ks-attractor-comparison}
\end{figure}

\paragraph{Qualitative attractor comparison.}
Figure~\ref{fig:ks-attractor-comparison} compares reference and learned
long-time trajectories of this same 14-node periodic system. Both start at
the final training snapshot and run for 1,000 additional time units, with the
first 200 discarded. A principal-component basis fitted only to the retained
reference states is used for both panels, with identical axis limits and
viewing angles. The plots compare sampled geometry after trajectory-wise
forecast agreement is lost. The panels occupy comparable coordinate ranges
but show visibly different sampled paths and occupancy; this figure alone
does not establish statistical agreement. They are numerical projections, not a proof of
chaos, equality of invariant measures, or a statement about the different
Neumann KS problem used for certification.

Using 149 common starts and fixing the training mean finds three matched
equilibria in each vector field (Table~\ref{tab:ksequilibria}). Field residual
and mean-error tolerances are $10^{-8}$; roots are deduplicated modulo discrete
shifts and matched within $10^{-3}$. Stable, unstable, and neutral counts agree
at real-part tolerance $10^{-5}$. These finite searches assess local fidelity
in the semidiscretization.

\begin{table}[htb]
\centering\small
\caption{Periodic KS equilibria. Position errors minimize over discrete shifts;
spectral errors use one-to-one eigenvalue assignment. Counts $(s,u,n)$ denote
stable, unstable, and neutral directions and agree for fitted and reference fields.}
\label{tab:ksequilibria}
\begin{tabular}{lrrc}
\toprule
Equilibrium & Position error & Spectral error & $(s,u,n)$\\
\midrule
Nearly constant & $1.02\times10^{-6}$ & $4.46\times10^{-6}$ & $(7,6,1)$\\
Nonconstant 1 & $4.12\times10^{-7}$ & $4.48\times10^{-6}$ & $(11,2,1)$\\
Nonconstant 2 & $8.34\times10^{-6}$ & $7.78\times10^{-5}$ & $(12,1,1)$\\
\bottomrule
\end{tabular}
\end{table}

Five surrogate-generated relative-orbit proposals are corrected in the reference
semidiscretization, with independent validation defects between approximately
$10^{-13}$ and $6\times10^{-12}$. Return up to spatial translation requires a
shift or drift variable and both temporal and spatial phase conditions.
These computations establish numerical closure in the 14-node system.

\subsection{Neumann-potential KS: orbit construction}
\label{sec:neumann-results}
On $[0,\ell]$, $\ell=22$, the mean-zero potential satisfies
\begin{align}
 h_t&=-\tfrac12h_x^2-h_{xx}-h_{xxxx}
       +\left\langle\tfrac12h_x^2\right\rangle_x,
 \label{eq:kspotential}\\
 h_x(0,t)&=h_x(\ell,t)=h_{xxx}(0,t)=h_{xxx}(\ell,t)=0.\nonumber
\end{align}
The velocity $u=h_x$ has $u=u_{xx}=0$ at the endpoints and an odd periodic
extension of length $44$. This differs from the length-22 periodic forecast
problem.

The 60-mode quadratic fit uses 7,324 training snapshots and 512 independent
derivative-evaluation snapshots, with relative error $5.98\times10^{-16}$
(rounded upward). This measures recovery of a quadratic Galerkin operator,
rather than time-series prediction. A recurrence from a 16-mode fitted
trajectory seeds continuation through 60 modes (Table~\ref{tab:continuation});
neither training nor initialization uses a known reference periodic orbit.

\begin{table}[htb]
\centering\small
\caption{Neumann KS continuation with fitted quadratic models. $M$ and $K$
are temporal and spatial orders; residuals have stage-specific scaling.}
\label{tab:continuation}
\begin{tabular}{rrrr}
\toprule
$M$ & $K$ & Physical period & Surrogate residual RMS\\
\midrule
16&16&73.36969486&$8.65\times10^{-17}$\\
24&20&72.70716150&$9.15\times10^{-17}$\\
40&28&72.70549994&$9.24\times10^{-17}$\\
64&36&72.70550498&$8.67\times10^{-17}$\\
80&40&72.70550498&$8.37\times10^{-12}$\\
96&44&72.70550498&$4.67\times10^{-12}$\\
128&52&72.70550498&$4.10\times10^{-12}$\\
160&60&72.70550498&$3.68\times10^{-15}$\\
\bottomrule
\end{tabular}
\end{table}

The final candidate has period approximately $72.705505$, independently
oversampled PDE RMS residual $2.42\times10^{-13}$, and closure error
$2.92\times10^{-11}$ using 72 modes and 40,000 steps per period.
These checks differ from the stage residuals and supply input to rigorous
validation, which additionally requires inverse, nonlinear, and infinite-tail
bounds.

Together, these experiments support finite-time forecasting and recovery of
selected local structures. Global fidelity depends on model selection, and
the Neumann orbit candidate supplies a separate input to rigorous validation;
forecast accuracy, long-time dynamics, and exact existence require distinct
evidence.

\paragraph{Numerical-assessment implementation.}
The forecast, attractor, Lyapunov, and equilibrium computations are implemented
in the Mackey--Glass and Kuramoto--Sivashinsky notebooks in the
\href{https://github.com/Aloh5/learned-dynamics-cap}{project's GitHub repository}.
The periodic KS orbit-transfer routines are provided in
\path{dpsr/ks_orbits.py}.

\section{Surrogate validation through certified invariant objects}
\label{sec:assessment}
We propose a method for validating learned surrogate models by testing whether
they generate invariant structures that can be certified in the reference
system. The method complements held-out prediction errors and numerical
attractor diagnostics with a local, mathematically controlled comparison.
Its contribution is the use of surrogate-generated candidates and their
certificates as evidence about the learned model's structural fidelity.
Polynomial fitting, invariant-solution searches, and a posteriori existence
theorems are established components \cite{brunton2016,page2021,framework2017}.
Here they are organized into a validation procedure that keeps the source of
the candidate and the scope of the resulting model assessment explicit.

\subsection{Validation procedure}
After fitting a vector field $\widehat g$ and evaluating its forecasts and
numerical dynamics, search for recurrence and refine a periodic candidate
using the surrogate. Export that candidate into the coordinates required by
the original equation, recording any compatibility constraints. Freeze the
waveform, frequency, coordinate conversion, and phase convention before the
proof stage. Then construct the original phase-augmented equation $F$ and
compute a posteriori bounds around this fixed export. Successful certification
gives
\begin{equation}
 \exists z_*:\quad F(z_*)=0,\qquad \|z_*-\bar z\|_X\leq r.
 \label{eq:object-assessment}
\end{equation}
The validation record specifies the candidate's provenance, the norm and
radius, fixed physical parameters, phase and symmetry restrictions, and
coverage of arithmetic errors and infinitely many omitted modes.

Freezing the export makes the certificate informative about the surrogate's
output: the enclosure is centered at the candidate the learned model supplied.
A correction in the reference equation could still produce a useful rigorous
existence result, but its distance from the original surrogate candidate would
also need to be tracked for this model-validation purpose. The primary
certificates therefore apply no such correction. Newton-refined certificates
are reported separately in Section~\ref{sec:cap-results}; there, interval
containment of the refined enclosure in the uncorrected one tracks the
required distance and identifies the same solution. This does not require zero access
to the governing equations: KS derivative labels and both proof operators use
the reference equations. The method therefore applies when an appropriate
reference operator and computable rigorous bounds are available.

\subsection{Interpretation alongside numerical validation}
A successful certificate validates one aspect of the surrogate: its ability
to supply a specified recurrent structure close to an exact object of the
reference system. The radius bounds the exported object's error in the stated
norm. It does not bound global vector-field error or indefinitely synchronized
trajectory error. In particular, the delay-consistent Mackey--Glass export
need not be an exact orbit of the unconstrained surrogate; the certificate
assesses the candidate produced by the stated construction procedure.

This interpretation preserves the role of the other results. Forecasts test
trajectory agreement over finite intervals, while numerical attractor
projections, distributions, Lyapunov estimates, and equilibrium spectra examine
selected dynamical properties. Certification adds a rigorous local test; it
does not replace those assessments or establish global equivalence of flows.
A model can therefore supply a certifiable periodic candidate while failing
long-time rollout tests, as the degree-five Mackey--Glass example illustrates.

Proof radii depend on scaling, inverse conditioning, and the sharpness of the
bounds, so they are not a universal ranking of models across different norms
or systems. A certification fraction would additionally require a candidate
sampling rule, duplicate handling, and common search and proof budgets. Here
we report two successful applications, not a population success rate or a
controlled comparison of validation methods. Failure of a sufficient
contraction criterion is inconclusive about whether a nearby exact object
exists. These qualifications define what the method validates without
weakening the quantitative statement in \eqref{eq:object-assessment}.

\section{Candidate construction and certification bounds}
\label{sec:cap}

\subsection{Harmonic balance and candidate export}
For a fitted ODE $\dot X=\widehat g(X)$, represent a periodic candidate as
\begin{equation}
 X_N(\theta)=b_0+\sum_{n=1}^{N}
 \bigl(b_n\cos(n\theta)+c_n\sin(n\theta)\bigr),\qquad \theta=\omega t.
\end{equation}
A close return seeds the solve. Harmonic balance equates Fourier coefficients
of $\omega X_N'$ and $\widehat g(X_N)$; collocation equates grid values.
A temporal phase condition removes time-translation freedom. The waveform
and frequency are refined by Newton or least-squares steps in the fitted
model, with enough temporal resolution to avoid nonlinear aliasing.
Increasing the harmonic order does not change the polynomial degree of
$\widehat g$. Numerical convergence alone supplies no existence theorem;
the frozen export is validated independently in the original equation.

For Mackey--Glass, a further compatibility issue arises. The six components of a free surrogate orbit need not be delayed copies of a single scalar history. The final export therefore parameterizes one scalar trigonometric polynomial $w$ and reconstructs all components from it:
\begin{equation}
 X_j(\theta)=w\!\left(\theta-(j-1)2\pi\Delta/T\right),
 \qquad j=1,\ldots,6.
 \label{eq:delaycompatible}
\end{equation}
The degree-five fitted ODE is evaluated on these constrained states. The least-squares refinement uses 40 temporal harmonics, 384 phase samples, a temporal phase condition, and penalties preserving the seed's mean and oscillation energy. The period is restricted to within five percent of the seed. Since these conditions need not admit an exact fitted-flow orbit, the exported waveform is described as a delay-consistent periodic candidate. Its reported surrogate collocation RMS is $0.0038437$, and its physical period is $46.7310078183329$. The construction audit records zero pre-CAP calls to the true DDE.

For Neumann KS, all discovery and continuation optimizations evaluate fitted vector fields. Reality and the observed half-period reflection symmetry are canonicalized at export, changing coefficients by at most $6.47\times10^{-14}$. The candidate is then frozen. ``Zero true-system Newton corrections'' means that no Newton update of this waveform is performed using the true equation; it does not mean that no surrogate optimization, Fourier conversion, or phase choice occurred. Nor does it mean zero access to the original equations: the KS training labels are exact spectral derivatives, and the frozen candidate is evaluated in the true equation for diagnostics and proof.

\subsection{An a posteriori contraction criterion}
Write the phase-augmented, suitably smoothed invariance equation as $F(z)=0$
on a Banach space $X$. A bounded injective approximate inverse $A$ gives the
Newton-like map $\mathcal T(z)=z-AF(z)$. We use the following scalar form of
the radii-polynomial (Newton--Kantorovich-type) contraction theorem; the
Fourier/radii-polynomial methodology is developed in
\citet{lessard2016,framework2017,berg2022}. We state the criterion and the
bounds used, rather than repeat the general derivations.

\begin{proposition}[Scalar-radius validation]
\label{prop:radii}
Suppose $F$ is continuously differentiable on a neighborhood of
$B_{r_*}(\bar z)$ and $A$ is bounded and injective. Assume
\begin{align}
 \|AF(\bar z)\|_X&\leq Y,\\
 \|I-ADF(\bar z)\|&\leq\alpha,\\
 \|A[DF(\bar z+v)-DF(\bar z)]\|&\leq2Z_2\|v\|_X
 \quad(\|v\|_X\leq r_*).
\end{align}
If $0<r\leq r_*$ satisfies
\begin{equation}
 p(r)=Y-(1-\alpha)r+Z_2r^2<0,
 \qquad \alpha+2Z_2r<1,
 \label{eq:scalar-radii}
\end{equation}
then $F$ has a unique zero in $B_r(\bar z)$.
\end{proposition}
Here $Y$ bounds the preconditioned residual and
$Z(r)=\alpha+2Z_2r$ bounds the derivative of $\mathcal T$ throughout the ball.
Some references write $Z_1$ for $\alpha$ and use a different normalization of
the quadratic constant; the factor of two in $2Z_2$ matters when comparing
certificates. Finite inverse defects, all finite/tail couplings, and the
infinite remainder are included. Uniqueness is within the phase slice and any
imposed symmetry space. These references supply the general theorem, not the
numerical bounds for our particular candidates.

\subsection{Mackey--Glass with fixed physical parameters}
\label{sec:mgcap}
We validate the candidate directly in \eqref{eq:mgphysical}, fixing the feedback coefficient, damping, delay, and exponent. Set $s=t/10$, $\tau=17/10$, and $G(z)=z/(1+z^{10})$. With $\theta=\omega s$ and $w(\theta)=\sum_{k\in\mathbb Z}u_ke^{ik\theta}$, the periodic equation is
\begin{equation}
 \omega w'(\theta)+w(\theta)-2G(w(\theta-\omega\tau))=0.
 \label{eq:mgproof}
\end{equation}
The unknown scalar is $\omega$; the physical period is $T=20\pi/\omega$. In particular, no damping parameter is varied. We work in the real Banach space of conjugate-symmetric Fourier coefficients with
\begin{equation}
 \|u\|_\nu=\sum_{k\in\mathbb Z}|u_k|\nu^{|k|},\quad
 \nu=1.01,\qquad
 \|(\delta\omega,h)\|_X=|\delta\omega|+\|h\|_\nu.
 \label{eq:mgnorm}
\end{equation}
Fourier validation in weighted sequence spaces follows the general approach of \citet{berg2022}; the rational nonlinearity is treated here by verified partial-fraction reciprocals.

Let $\omega_{\rm ref}>0$ be the candidate frequency, held constant as an equation-scaling factor. Let $P$ remove the zero Fourier mode, define $(S_\omega a)_k=e^{-ik\omega\tau}a_k$, and set $J_0=1$, $J_k=(ik\omega_{\rm ref})^{-1}$ for $k\neq0$. The integrated state residual is
\begin{equation}
 F_u(\omega,u)=\frac{\omega}{\omega_{\rm ref}}Pu
                  +J\bigl(u-2S_\omega G(u)\bigr).
 \label{eq:mgfixedoperator}
\end{equation}
Its zero component is $u_0-2G(u)_0$; each nonzero component is \eqref{eq:mgproof} divided by $ik\omega_{\rm ref}$. A finite linear phase condition $\varphi u-b=0$ completes $F$. The row $\varphi$ is obtained from the conjugate temporal tangent of the candidate, normalized in the finite Euclidean inner product, with support $|k|\leq80$. Both $\varphi$ and $b$ are stored as dyadic data, and their rounding residual is included in the proof. This phase condition preserves the time origin of the exported waveform.

The approximate inverse acts on the augmented core $|k|\leq80$ and equals
the identity on the state tail. The finite inverse is verified by an interval
defect bound below $6.71\times10^{-12}$. The implementation includes the
frequency derivative of the delayed phase; fixing the physical delay does not
fix the phase shift $\omega\tau$.

For the rational nonlinearity, verified degree-160 Fourier reciprocals of
$\bar w-p_j$, where $p_j=e^{i\pi(2j-1)/10}$, enclose the partial fractions of
$G$. Let $B_j$ bound the corresponding exact reciprocal norm. The checks
$B_jr<1$ exclude every pole on the full proof ball, by a Neumann-series
argument for the reciprocal in the Banach algebra of Fourier coefficients. All
physical parameters are fixed, and the frequency is the extra unknown.

The actual estimates are $Y\leq4.194761221\times10^{-4}$,
$\alpha=Z_0\leq0.069980673$, and $2Z_2=C\leq494.208949$ on the
radius-$0.0006$ ball. The full infinite-tail contribution to the linear defect
is bounded by $0.068703474$, and the maximum $B_jr$ is below $0.003483453$.
Table~\ref{tab:mgcap} and Theorem~\ref{thm:mgorbit} record the resulting
inequalities. The interval convolution and weighted-norm calculations,
including the frequency column and omitted modes, are implemented in
\path{CAP/Mackey-Glass_PO/proof_MG_fixed_parameter.m}.

\subsection{Neumann KS: scaled space--time Fourier equation}
\label{sec:kscap}
Set
\begin{equation}
 q=\pi/22,\qquad y=qx,\qquad s=q^2t,\qquad v=u/q,
 \qquad \nu=q^2.
\end{equation}
Then $v_s+vv_y+v_{yy}+\nu v_{yyyy}=0$ on $0\leq y\leq\pi$. With temporal phase $\theta=\omega t$ and scaled frequency $f=\omega/q^2$, the periodic equation is
\begin{equation}
 E(f,v)=f v_\theta+\nu v_{yyyy}+v_{yy}+\tfrac12\partial_y(v^2)=0.
 \label{eq:ksfourier}
\end{equation}
The unknowns are $f$ and the space--time Fourier coefficients of $v$. A real basis uses a constant temporal term and cosine/sine temporal amplitudes multiplying spatial sine modes. Only the temporal phase is fixed: a generic spatial translate violates the endpoint conditions. The candidate satisfies
\begin{equation}
 v(\theta+\pi,y)=-v(\theta,\pi-y),
\end{equation}
so the proof is carried out in the closed invariant coefficient subspace with $n+k$ even, where $n$ and $k$ index temporal and spatial modes.

If $a_{0k},a^c_{nk},a^s_{nk}$ are the real amplitudes, the correction norm is
\begin{equation}
 \|(\delta f,\delta v)\|_X=|\delta f|
 +\sum_k|\delta a_{0k}|
 +\sum_{n\geq1,k\geq1}(|\delta a^c_{nk}|+|\delta a^s_{nk}|),
 \label{eq:ksnorm}
\end{equation}
restricted to the allowed parity. This is an unweighted real trigonometric coefficient $\ell^1$ norm, not a coefficient supremum norm and not the weighted complex norm of \eqref{eq:mgnorm}. Product-to-sum identities bound products in the corresponding trigonometric coefficient algebra, since the absolute coefficients of each elementary product sum to at most one. The spatial and temporal derivative operators are handled after smoothing.

At the candidate $(\bar f,\bar v)$ define
\begin{equation}
 L_0=\bar f\partial_\theta+\nu\partial_y^4+\partial_y^2,
 \quad c=1/\nu,\quad S_c=L_0+c.
\end{equation}
Define the zero-finding operator explicitly by
\begin{equation}
 F_{\rm KS}(f,v)=\begin{pmatrix}
 \varphi(v-\bar v)\\ S_c^{-1}E(f,v)
 \end{pmatrix},
 \qquad T=\frac{2\pi}{q^2f}.
 \label{eq:kszero}
\end{equation}
The finite linear functional $\varphi$ uses the candidate's temporal tangent.
The smoothing operator $S_c$ is fixed at the candidate. Thus
$F_{\rm KS}=0$ enforces both the PDE and the chosen time origin. The frequency
unknown balances this added scalar phase equation. Reality, odd spatial
extension, and parity are enforced in the coefficient space.

Following \citet{figueras2017}, the derivative at the candidate is
\begin{equation}
 DF_{\rm KS}(\bar f,\bar v)=
 \begin{pmatrix}
 0&\varphi\\
 S_c^{-1}\bar v_\theta&I-cS_c^{-1}+S_c^{-1}\partial_y(\bar v\,\cdot)
 \end{pmatrix}.
 \label{eq:kslinearization}
\end{equation}
Project this operator onto the frequency/phase variables and modes
$0\leq n\leq160$, $1\leq k\leq40$, and numerically invert the resulting
6,421-dimensional matrix to obtain $B$. Extend $B$ to an approximate inverse
$A$ by using $L_0^{-1}S_c$ on the complementary state modes. Equivalently,
the tail acts on the unsmoothed PDE residual with $L_0^{-1}$, whose complex
Fourier symbol is $(in\bar f+\nu k^4-k^2)^{-1}$. This explains how a finite
matrix inverse becomes a Newton-like inverse for the infinite equation.
The inverse core has 40 spatial modes, although the frozen candidate has 60.

The map $(f,v)\mapsto(f,v)-AF_{\rm KS}(f,v)$ is used to prove contraction
near the frozen candidate, without applying a true-equation Newton update to
that candidate. The retained nonlinear correction is exactly
$\delta f\,\partial_\theta\delta v+\tfrac12\partial_y(\delta v)^2$.
Preconditioning and coefficient-algebra estimates give the certified bounds
\begin{equation}
 Y\leq4.119405\times10^{-7},\qquad
 \alpha\leq0.890809345,\qquad Z_2\leq165.085084.
 \label{eq:ks-method-bounds}
\end{equation}
These are candidate-specific computed estimates, not constants taken from
the reference paper.

\paragraph{Coverage and arithmetic.}
The quadratic residual is enclosed using 641 temporal nodes and 120 spatial
output modes. The linear bound includes the finite frequency column and every
finite/tail coupling: the complete halo through $n=320$, $k=100$, the remaining
temporal gap through $n=359$, and analytic symbol bounds for all larger
indices. Thus the estimate includes infinitely many modes, not just a larger
Galerkin truncation. The derivative bounds likewise include their infinite
tails. The stored finite inverse is checked through its defect, and the
preconditioner is injective on the phase-fixed space.

INTLAB encloses the residual and final scalar inequalities \cite{rump1999}.
The large finite matrices and coupling scans use binary64 arithmetic with
explicit forward-error majorants \cite{higham2002}, as detailed in
\path{CAP/Kuramoto_PO/PROOF_AUDIT.md}. Those guards depend on the audited
truncations and arithmetic assumptions. The driver checks that saved scan data
match the frozen candidate. General Fourier and smoothing estimates are
referenced to \citet{figueras2017,framework2017}; the implementation and audit
supply the adaptations to this boundary condition, symmetry space, norm, and
candidate. No claim is made that a citation alone validates those adaptations.

\section{Certified periodic-orbit results}
\label{sec:cap-results}
\begin{figure}[htb]
\centering
\includegraphics[width=\textwidth]{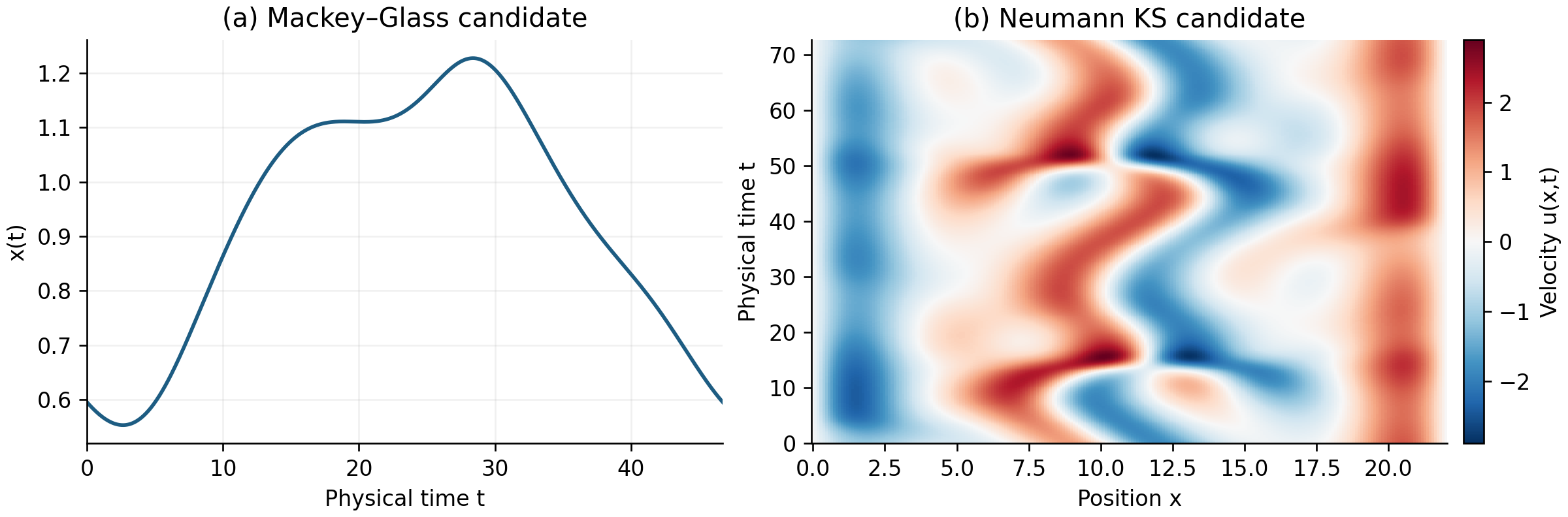}
\caption{Frozen surrogate-generated candidates used in the certificates.
Left: the Mackey--Glass scalar waveform over one candidate period.
Right: physical Neumann KS velocity $u(x,t)$ over one candidate period.
These are plots of the exported approximations, not exact solution curves;
Theorems~\ref{thm:mgorbit} and \ref{thm:ksorbit} enclose nearby exact solutions
in their respective augmented coefficient norms.}
\label{fig:certified-candidates}
\end{figure}

The following results apply the validation procedure of
Section~\ref{sec:assessment}, using the operators and bounds in
Section~\ref{sec:cap}. They quantify the accuracy of two frozen,
surrogate-generated periodic candidates at fixed physical parameters and
provide tighter enclosures after Newton refinement in the original
equations.
The stated phase and symmetry restrictions specify the local uniqueness
claims. Together with the forecasting and qualitative results, these
certificates provide complementary evidence about the learned models.
Figure~\ref{fig:certified-candidates} displays both candidates.

\subsection{Fixed-parameter Mackey--Glass periodic orbit}
\begin{samepage}
\begin{theorem}[Computer-assisted Mackey--Glass periodic orbit]
\label{thm:mgorbit}
The fixed-parameter equation \eqref{eq:mgphysical} has a nonstationary periodic solution. Relative to the uncorrected surrogate candidate and the phase condition above, its scaled frequency and Fourier coefficients satisfy
\begin{equation}
 |\omega_*-\bar\omega|+\|u_*-\bar u\|_{1.01}
 \leq6.00\times10^{-4}.
 \label{eq:mgcertifiedradius}
\end{equation}
It is the unique solution in this ball and phase slice. Its physical period satisfies
\begin{equation}
 T\in[46.7101634894597,\;46.7518707590022].
 \label{eq:mgperiodinterval}
\end{equation}
\end{theorem}
\end{samepage}
\begin{proof}
INTLAB evaluation gives the bounds in Table~\ref{tab:mgcap}. At $r=0.0006$, the largest pole-neighborhood factor $B_jr$ is below $0.003484$, so all nonlinear estimates hold on the full ball. The preconditioned residual and derivative estimates satisfy
\begin{align}
 Y+(Z_0-1)r+\tfrac12Cr^2&\leq-4.957786349\times10^{-5},\\
 Z_0+Cr&\leq0.366506043<1.
\end{align}
Proposition~\ref{prop:radii} applies with $\alpha=Z_0$ and $Z_2=C/2$. The validated finite inverse and identity tail make $A$ injective, hence the fixed point solves \eqref{eq:mgproof}. Conjugate symmetry gives a real solution. A nonconstant Fourier coefficient has modulus greater than $0.14895$ throughout the ball, excluding an equilibrium. Finally, $\bar\omega-r>0$ and outward-rounded evaluation of $20\pi/[\bar\omega-r,\bar\omega+r]$ gives \eqref{eq:mgperiodinterval}.
\end{proof}

\begin{table}[htb]
\centering
\caption{Fixed-parameter Mackey--Glass certificate. Upper bounds are rounded outward. The radius uses the augmented norm \eqref{eq:mgnorm}; its state component bounds the distance from the uncorrected surrogate waveform.}
\label{tab:mgcap}
\begin{tabular}{lr}
\toprule
Quantity & Value or upper bound\\
\midrule
True-system Newton corrections &0\\
Candidate physical period &$46.7310078183329$\\
Preconditioned residual $Y$ &$4.194761221\times10^{-4}$\\
Linear defect $Z_0$ &$0.069980673$\\
Derivative-variation constant $C$ &$494.208949$\\
Maximum reciprocal defect &$3.90224\times10^{-14}$\\
Maximum pole-neighborhood factor $B_jr$ &$0.003483453$\\
Certified augmented radius &$6.00\times10^{-4}$\\
Contraction constant &$0.366506043$\\
\bottomrule
\end{tabular}
\end{table}

The proof uses direct outward-rounded convolution, matrix products, weighted norms, and scalar inequalities in MATLAB with INTLAB \cite{rump1999}. Floating-point Fourier calculations supply only approximate reciprocals and the finite inverse, each verified by its interval defect. The certificate therefore establishes fixed-parameter existence directly from the surrogate export with zero true-equation Newton updates. The enclosed period need not be minimal, and the theorem makes no assertion about orbital stability or chaotic dynamics.

\subsection{Newton-refined Mackey--Glass certificate}
The uncorrected enclosure measures the accuracy of the surrogate-generated
candidate. To resolve the same periodic solution more precisely, we also
apply five Newton updates in the original phase-augmented Mackey--Glass
Fourier equation. Physical parameters, the reference frequency used for
equation scaling, and the phase condition remain fixed. Only the candidate
frequency and waveform are corrected. The numerical inverse uses order 80,
and reciprocal and derivative approximations use order 160, as in the
uncorrected proof. All interval bounds are recomputed at the corrected center.
On complementary state modes the approximate inverse multiplies by
$\omega_{\rm ref}/\widetilde\omega$, accounting for the corrected frequency
in the integrated operator.

\begin{theorem}[Newton-refined Mackey--Glass periodic orbit]
\label{thm:mgorbit-refined}
Let $(\widetilde\omega,\widetilde u)$ be the candidate obtained by five
Newton updates in the phase-augmented, fixed-parameter Mackey--Glass
Fourier equation, starting from the candidate of Theorem~\ref{thm:mgorbit}.
The nonstationary periodic solution $(\omega_*,u_*)$ of
Theorem~\ref{thm:mgorbit} satisfies
\[
 |\omega_*-\widetilde\omega|+
 \|u_*-\widetilde u\|_{1.01}\leq 1.2\times10^{-12}.
\]
It is unique in this ball and the stated phase slice. Its physical period satisfies
\begin{equation}
 T\in[46.72953965885938,\;46.72953965894285].
 \label{eq:mgrefinedperiodinterval}
\end{equation}
\end{theorem}
\begin{proof}
Outward-rounded INTLAB calculations give the bounds in
Table~\ref{tab:mgcap-refined}.
At $r=1.2\times10^{-12}$, the validated radii polynomial is negative
and the contraction bound is less than one. The proof includes all Fourier
tails, pole exclusion, nonstationarity, and injectivity of the approximate
inverse, so Proposition~\ref{prop:radii} applies.
The phase condition is identical to that of the uncorrected certificate.
Interval arithmetic bounds the distance from the original center to the
entire refined ball by $4.194907107\times10^{-4}
<6.00\times10^{-4}$. Uniqueness in the original ball identifies the two
certified solutions. The period enclosure follows from outward-rounded
frequency conversion.
\end{proof}

\begin{table}[htb]
\centering
\caption{Newton-refined Mackey--Glass certificate in the augmented norm
\eqref{eq:mgnorm}. Proof bounds are rounded outward. The numerical residual
is an unweighted finite-core residual and is not the certified distance.}
\label{tab:mgcap-refined}
\begin{tabular}{lr}
\toprule
Quantity & Value or upper bound\\
\midrule
True-system Newton corrections &5\\
Corrected candidate physical period &$46.7295396589011$\\
Numerical finite-core residual $\ell^1$ norm &$9.15\times10^{-17}$\\
Preconditioned residual $Y$ &$1.0332785\times10^{-12}$\\
Linear defect $Z_0$ &$0.070006496$\\
Derivative-variation constant $C$ &$489.448942$\\
Certified augmented radius &$1.2\times10^{-12}$\\
Contraction constant &$0.0700065$\\
Distance from uncorrected center to refined ball &$4.194907107\times10^{-4}$\\
\bottomrule
\end{tabular}
\end{table}

The certified radius is $5\times10^8$ times smaller than the uncorrected
radius. These enclosures answer different questions: the first bounds the
error of the surrogate-generated export, while the second resolves the exact
orbit after correction in the original equation. Combined, they also sharpen
the first answer: since the exact solution lies in the refined ball, the
surrogate export is within $4.194907107\times10^{-4}$ of it in
\eqref{eq:mgnorm}, so the uncorrected radius $6.00\times10^{-4}$ is within
a factor of $1.5$ of the true export error. The numerical residual
$9.15\times10^{-17}$ is smaller than the certified distance because the proof
also encloses reciprocal approximation, arithmetic, and infinite-tail
contributions. Neither enclosure establishes orbital stability, minimal
period, or chaos.

\subsection{KS certificate and interpretation}
The main validation result concerns the candidate $(\bar f,\bar v)$ constructed by fitted-model continuation in Table~\ref{tab:continuation}.

\begin{samepage}
\begin{theorem}[Computer-assisted existence of a Neumann KS periodic solution]
\label{thm:ksorbit}
For $\ell=22$, the KS velocity equation with endpoint conditions $u=u_{xx}=0$ has a nonstationary time-periodic solution. In the scaled coordinates of Section~\ref{sec:kscap}, a representative $(f_*,v_*)$ satisfies
\begin{equation}
 \|(f_*-\bar f,v_*-\bar v)\|_X\leq10^{-5},
\end{equation}
where the norm is defined in \eqref{eq:ksnorm}. This solution is locally unique in that ball within the phase-fixed, half-period reflection-symmetric space. Its physical period belongs to the interval in \eqref{eq:ksperiodinterval}.
\end{theorem}
\end{samepage}
\begin{proof}
The computer-assisted estimates give the conservative upper bounds
\begin{equation}
 Y\leq4.119405\times10^{-7},\qquad
 \alpha\leq0.890809345,\qquad
 Z_2\leq165.085084.
 \label{eq:ksbounds}
\end{equation}
At $r=10^{-5}$, the outward-rounded computation records
\begin{equation}
 p(r)\leq-6.634575\times10^{-7},\qquad
 \alpha+2Z_2r\leq0.894111046<1.
\end{equation}
The chosen radius lies strictly between the recorded polynomial roots, approximately $3.79444\times10^{-6}$ and $6.57626\times10^{-4}$. Together with the finite and infinite estimates, these bounds satisfy Proposition~\ref{prop:radii} for the phase-fixed symmetry subspace. A nonzero time-dependent coefficient, bounded below by more than $4.16$, excludes an equilibrium within the radius. The frequency enclosure is positive and yields the outward-rounded physical-period interval
\begin{equation}
 T\in[72.7053334243950,\;72.7056765389189].
 \label{eq:ksperiodinterval}
\end{equation}
\end{proof}
The coefficient correction is bounded by $10^{-5}$ in \eqref{eq:ksnorm}; this also bounds the state part of that norm. Restoring physical velocity multiplies a uniform state bound by $q$, but the period is also enclosed, so this should not be interpreted as an indefinitely synchronized physical-time trajectory error.

The certificate reports complete tail validation, zero true-system optimizer calls during candidate construction, and zero true-system Newton corrections. It proves a nearby ordinary periodic solution of the Neumann boundary-value problem in the specified symmetry class. It does not certify orbital stability, a Floquet spectrum, chaotic dynamics, or a global equivalence of the fitted flow and the PDE. 

\subsection{Newton-refined Neumann KS certificate}
Starting from the uncorrected Neumann KS export, we extend the candidate to
220 temporal and 100 spatial modes and perform three Newton--GMRES updates
in the original Fourier equation. The physical interval, parameters,
half-period reflection symmetry, and temporal phase condition remain fixed.
The original finite inverse with orders $(160,40)$ and the original-frequency
tail inverse remain the preconditioner. The smoothing operator also remains
fixed. Thus, unlike the Mackey--Glass refinement, in which all bounds are
recomputed, the original full-space derivative estimates are transferred
to the corrected center without asserting a new finite/tail scan.

Let $z_0=(\bar f,\bar v)$ be the original center and
$\widetilde z=(\widetilde f,\widetilde v)$ the corrected candidate. For the
fixed preconditioned operator $H=AF_{\rm KS}$, the global quadratic remainder
bound gives
\begin{equation}
 \|I-DH(\widetilde z)\|\leq
 \alpha_0+2Z_{2,0}\|\widetilde z-z_0\|_X,
 \qquad Z_{2,\rm refined}=Z_{2,0}.
 \label{eq:ksrefined-transfer}
\end{equation}
The phase row is unchanged. The distance between centers and the new
preconditioned residual are evaluated with interval arithmetic. Direct
coefficient convolutions and exact spatial sine product identities enclose
the complete quadratic residual through temporal order 440 and spatial order
200. This avoids the trigonometric-transform residual enclosure used for the
uncorrected certificate. All higher residual coefficients at the finite
center vanish exactly, while the transferred derivative bounds still cover
infinitely many correction modes.

\begin{theorem}[Newton-refined Neumann KS periodic solution]
\label{thm:ksorbit-refined}
Let $(\widetilde f,\widetilde v)$ be the candidate obtained by three
Newton--GMRES updates in the phase-augmented Neumann KS Fourier equation,
starting from the candidate of Theorem~\ref{thm:ksorbit}. The nonstationary
periodic solution $(f_*,v_*)$ of Theorem~\ref{thm:ksorbit} satisfies
\begin{equation}
 \|(f_*-\widetilde f,v_*-\widetilde v)\|_X\leq7\times10^{-9},
 \label{eq:ksrefinedradius}
\end{equation}
where $\|\cdot\|_X$ is the real trigonometric coefficient norm
\eqref{eq:ksnorm}. It is unique in this ball within the same phase-fixed,
half-period reflection-symmetric space. Its physical period satisfies
\begin{equation}
 T\in[72.7055048611696,\;72.7055051013500].
 \label{eq:ksrefinedperiodinterval}
\end{equation}
\end{theorem}
\begin{proof}
INTLAB encloses the center displacement by
$\|\widetilde z-z_0\|_X\leq4.915052\times10^{-10}$.
Using \eqref{eq:ksrefined-transfer} and the direct interval residual gives
the bounds in Table~\ref{tab:kscap-refined}. At $r=7\times10^{-9}$,
outward-rounded evaluation yields
\begin{align}
 Y-(1-\alpha)r+Z_2r^2&\leq-6.7893545\times10^{-11},\\
 \alpha+2Z_2r&\leq0.890811818<1.
\end{align}
Proposition~\ref{prop:radii} applies with the original injective
preconditioner. A time-dependent coefficient has modulus greater than
$4.1674$ throughout the ball, excluding an equilibrium. The full refined
ball lies within distance $7.491506\times10^{-9}<10^{-5}$ of the original
center. Since the equation and phase/symmetry slice agree, uniqueness in the
original ball identifies the same periodic solution. Positivity of the
frequency enclosure and outward-rounded conversion give
\eqref{eq:ksrefinedperiodinterval}.
\end{proof}

\begin{table}[htb]
\centering
\caption{Newton-refined Neumann KS certificate in the norm
\eqref{eq:ksnorm}. Proof bounds are rounded outward. The numerical residual
is the finite preconditioned residual used during refinement, not a distance
bound.}
\label{tab:kscap-refined}
\begin{tabular}{lr}
\toprule
Quantity & Value or upper bound\\
\midrule
True-system Newton corrections &3\\
Candidate temporal/spatial orders &220 / 100\\
Finite inverse temporal/spatial orders &160 / 40\\
Corrected candidate physical period &$72.7055049812598$\\
Numerical preconditioned core residual &$1.88\times10^{-13}$\\
Certified preconditioned residual $Y$ &$6.964318202\times10^{-10}$\\
Transferred derivative defect $\alpha$ &$0.890809507$\\
Quadratic bound $Z_2$ &$165.085083326$\\
Certified augmented radius &$7\times10^{-9}$\\
Contraction constant &$0.890811818$\\
Distance from uncorrected center to refined ball &$7.491506\times10^{-9}$\\
\bottomrule
\end{tabular}
\end{table}

The certified radius is approximately 1,400 times smaller than the
uncorrected $10^{-5}$ enclosure. Newton refinement, the increased candidate
orders, and the sharper residual enclosure all contribute to the result;
these computations do not isolate their individual effects. The numerical
residual is much smaller than the certified distance because the latter
also accounts for conditioning, interval arithmetic, and all Fourier tails.
The original certificate continues to quantify the accuracy of the
uncorrected surrogate export. The refined certificate sharpens that
quantification: containment places the exact solution within
$7.491506\times10^{-9}$ of the uncorrected export in \eqref{eq:ksnorm}. The
uncorrected radius $10^{-5}$ is therefore conservative by roughly three
orders of magnitude, consistent with the conservative trigonometric-transform
residual bound and the chosen proof radius of that certificate.

\subsection{Neumann KS from station time series}
\label{sec:station-results}
\begin{figure}[htb]
\centering
\includegraphics[width=\textwidth]{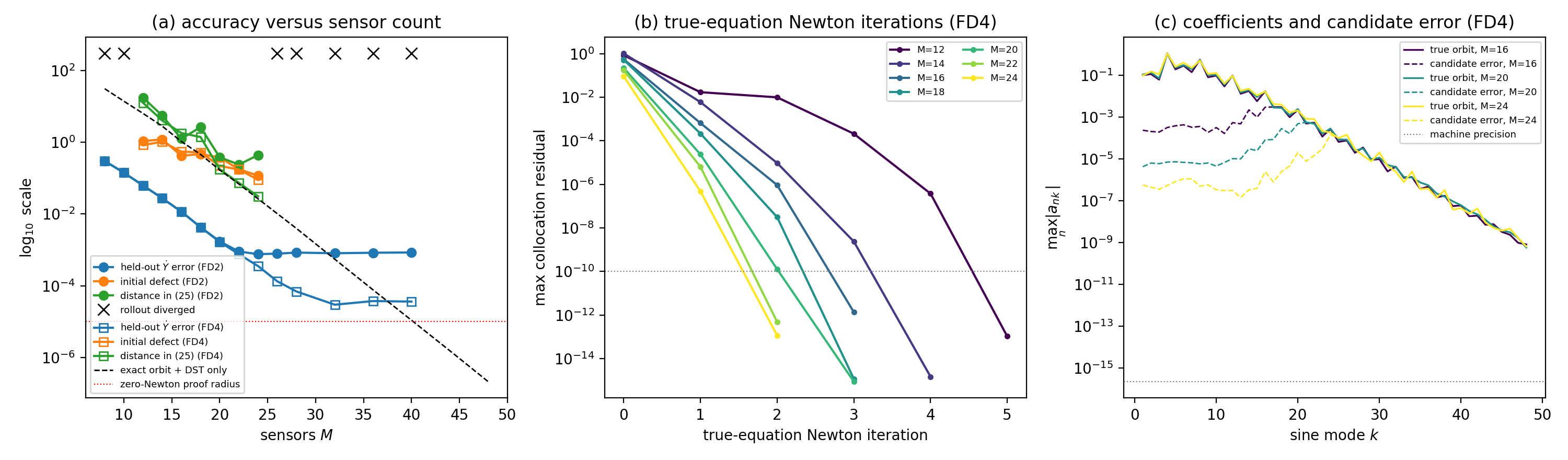}
\caption{Neumann KS station surrogates. (a) Held-out derivative error, initial
defect, and distance to the Newton limit in \eqref{eq:ksnorm} versus sensor
count, with the distance obtained from the exact orbit and the DST alone
(dashed). Crosses mark sensor counts at which every rollout diverged. The
dotted line is the radius $10^{-5}$ of the uncorrected certificate.
(b) Maximum true collocation residual during Newton iteration from each frozen
FD4 candidate. (c) Spatial coefficient magnitudes $\max_n|a_{nk}|$ of the
Newton limit (solid) and of the candidate error (dashed), for FD4 candidates
with $M=16,20,24$.}
\label{fig:station-convergence}
\end{figure}

We now apply the procedure to the station surrogates of
Section~\ref{sec:stations}, for $8\leq M\leq40$.
Table~\ref{tab:stations} and Figure~\ref{fig:station-convergence} summarize
the results. Most surrogates have one to three periodic orbits, and all 24
candidates converge under Newton iteration to true periodic orbits, with
periods $69.25$, $71.11$, $72.71$, $73.12$, $73.33$, $73.58$, $73.59$, $75.72$,
$85.68$ and $87.09$. For each $M$ the table reports the candidate whose Newton
limit has period $75.72$, the orbit certified in
Theorem~\ref{thm:ksorbit-stations} below. If there is none, it reports the
candidate closest to its own Newton limit. For example, with FD4 at $M=22$ and
$M=24$ no candidate converges to the period-$75.72$ orbit; the reported
candidates converge instead to the period-$72.71$ orbit of
Theorem~\ref{thm:ksorbit}, at the resolution of the test.

\begin{table}[htb]
\centering\small
\caption{Neumann KS station surrogates on the grid \eqref{eq:station-sites}.
$e$: held-out relative derivative error of the selected model (smallest over
$\lambda$ if every rollout diverged). Defect: maximum true collocation residual
of the frozen candidate. $N$: true-equation Newton corrections to residual
$<10^{-10}$. $d$: distance in \eqref{eq:ksnorm} from the frozen candidate to its
Newton limit (32 temporal harmonics, 48 sine modes). ``div.'': every rollout
diverged. Last column: $d$ for the exact orbit sampled at the sites and
reconstructed by \eqref{eq:dst}.}
\label{tab:stations}
\begin{tabular}{r@{\quad}rrrr@{\quad}rrrr@{\quad}r}
\toprule
 & \multicolumn{4}{c}{FD2 derivatives} & \multicolumn{4}{c}{FD4 derivatives} & DST only\\
$M$ & $e$ & Defect & $N$ & $d$ & $e$ & Defect & $N$ & $d$ & $d$\\
\midrule
8  &$3.0\times10^{-1}$&\multicolumn{3}{c}{div.}&$3.0\times10^{-1}$&\multicolumn{3}{c}{div.}&$3.1\times10^{1}$\\
10 &$1.4\times10^{-1}$&\multicolumn{3}{c}{div.}&$1.4\times10^{-1}$&\multicolumn{3}{c}{div.}&$1.4\times10^{1}$\\
12 &$6.2\times10^{-2}$&1.08&9&17.6&$6.3\times10^{-2}$&0.84&5&12.6&6.36\\
14 &$2.7\times10^{-2}$&1.21&4&5.51&$2.8\times10^{-2}$&1.03&4&4.04&2.82\\
16 &$1.2\times10^{-2}$&0.42&3&1.25&$1.2\times10^{-2}$&0.55&3&1.80&1.05\\
18 &$4.3\times10^{-3}$&0.47&4&2.63&$4.3\times10^{-3}$&0.52&3&1.40&0.46\\
20 &$1.7\times10^{-3}$&0.38&3&0.380&$1.7\times10^{-3}$&0.22&3&0.175&0.17\\
22 &$9.0\times10^{-4}$&0.18&2&0.239&$7.6\times10^{-4}$&0.17&2&0.074&$6.8\times10^{-2}$\\
24 &$8.2\times10^{-4}$&0.12&3&0.438&$4.0\times10^{-4}$&0.09&2&0.031&$2.7\times10^{-2}$\\
26 &$7.9\times10^{-4}$&\multicolumn{3}{c}{div.}&$1.3\times10^{-4}$&\multicolumn{3}{c}{div.}&$1.0\times10^{-2}$\\
28 &$8.4\times10^{-4}$&\multicolumn{3}{c}{div.}&$7.0\times10^{-5}$&\multicolumn{3}{c}{div.}&$4.0\times10^{-3}$\\
32 &$8.1\times10^{-4}$&\multicolumn{3}{c}{div.}&$3.0\times10^{-5}$&\multicolumn{3}{c}{div.}&$5.5\times10^{-4}$\\
36 &$8.3\times10^{-4}$&\multicolumn{3}{c}{div.}&$3.8\times10^{-5}$&\multicolumn{3}{c}{div.}&$7.7\times10^{-5}$\\
40 &$8.5\times10^{-4}$&\multicolumn{3}{c}{div.}&$3.6\times10^{-5}$&\multicolumn{3}{c}{div.}&$1.1\times10^{-5}$\\
\bottomrule
\end{tabular}
\end{table}

\paragraph{Convergence in the number of sensors.}
For $12\leq M\leq24$ every surrogate is stable and yields periodic candidates.
Each candidate's period is within $1.3\%$ of its Newton limit's. The held-out
derivative error decreases roughly geometrically, by about one decade per five
sensors, and the distance from each candidate to its own Newton limit follows
it, from $17.6$ at $M=12$ to $3.1\times10^{-2}$ at $M=24$ (FD4); these limits
are not all the same orbit. The number of Newton corrections falls
from nine to two. The initial defect decreases more slowly, because it is
dominated by the high, strongly damped modes absent from the candidate.

Figure~\ref{fig:station-convergence}(c) locates the error. At $M=24$ the
candidate error in the resolved modes $k\leq20$ is $10^{-7}$--$10^{-5}$, so the
learned dynamics of those modes are accurate. The remaining error is the
truncation itself: for $k>M$ the candidate is zero, and the error equals the true
coefficients, still about $10^{-4}$ at $k=24$. The last column of
Table~\ref{tab:stations} uses the period-$75.72$ orbit. Sampling instead the
period-$72.71$ orbit, which the FD4 candidates at $M=22$ and $24$ reach, and
reconstructing by \eqref{eq:dst} gives distances $7.4\times10^{-2}$ and
$3.2\times10^{-2}$, essentially equal to those candidates' distances
$7.4\times10^{-2}$ and $3.1\times10^{-2}$. In this range the surrogate is no
longer the main source of error; the number of resolved modes is.

\paragraph{Failure beyond 24 sensors.}
Reaching the uncorrected proof radius $10^{-5}$ requires about $M=40$ even with
exact data (Table~\ref{tab:stations}, last column). At every $M\geq26$, however,
the rollout diverges for every ridge parameter. This happens although the
held-out derivative error keeps decreasing with FD4, to $3\times10^{-5}$. With
FD2 it stalls at the finite-difference floor near $8\times10^{-4}$. The added
modes are tiny along the attractor ($|a_k|\lesssim10^{-3}$ for $k\geq25$) and
relax on time scales $1/((kq)^4-(kq)^2)$, with $q=\pi/22$ as in
Section~\ref{sec:kscap}, far below the sampling interval
($6.7\times10^{-3}$ at $k=25$, $9.7\times10^{-4}$ at $k=40$). The data therefore
contain almost no information about their damping. For an FD4 fit at $M=40$,
written in sine coordinates through \eqref{eq:dst}, the learned linear
coefficient of each mode matches the true value $(kq)^2-(kq)^4$ to three digits up to
$k=15$. It is $-52$ instead of $-58$ at $k=20$, and $+2.7$ instead of $-150$ at
$k=25$. Above $k=25$ it is close to zero, instead of reaching $-1032$ at $k=40$.
These errors barely affect the held-out derivative error, but without the strong
damping the rollouts diverge. Increasing the sampling rate would require steps
well below $10^{-3}$ for $k\approx40$.

Thus, increasing the number of sensors improves the surrogate until about
$M=24$. Beyond that, a quadratic regression on station data becomes unstable
before it resolves enough modes for an uncorrected certificate. In this setting
we did not find a proxy accurate enough for the radius $10^{-5}$ from sensor
data alone. Certification from station data therefore follows the
Newton-refined route.

\paragraph{Newton-refined station certificate.}
We certify one station-generated orbit after true-equation refinement. The
candidate is the FD2 surrogate orbit at $M=24$ ($\lambda=10^{-10}$, surrogate
period $75.7558$). It satisfies the half-period reflection symmetry to relative
defect $4.0\times10^{-4}$. We continue it with Newton in the true Galerkin
equation through temporal/spatial orders $(32,32)$, $(48,36)$, $(64,40)$,
$(96,44)$, $(128,52)$ and $(160,60)$. We then repeat the complete finite/tail
scan and the INTLAB validation of Section~\ref{sec:kscap} at the corrected
center, recomputing all bounds.

\begin{theorem}[Neumann KS periodic solution from station data]
\label{thm:ksorbit-stations}
For $\ell=22$, the KS velocity equation with $u=u_{xx}=0$ at the endpoints has a
nonstationary time-periodic solution within distance $10^{-5}$, in
\eqref{eq:ksnorm}, of the Newton-refined station candidate. It is locally unique
in that ball within the phase-fixed, half-period reflection-symmetric space,
and its physical period satisfies
\begin{equation}
 T\in[75.7237325983517,\;75.7241047934699].
\end{equation}
\end{theorem}
\begin{proof}
The validation gives $Y\leq1.500091\times10^{-7}$, $\alpha\leq0.8801236$, and
$Z_2\leq115.99072$. At $r=10^{-5}$, $p(r)\leq-1.037156\times10^{-6}$ and
$\alpha+2Z_2r\leq0.8824434<1$. The radius lies between the polynomial roots,
approximately $1.25\times10^{-6}$ and $1.03\times10^{-3}$. A time-dependent
coefficient exceeds $2.71$ in modulus throughout the ball, which excludes an
equilibrium. Proposition~\ref{prop:radii} applies.
\end{proof}

This periodic solution differs from that of Theorem~\ref{thm:ksorbit}. The frozen
station candidate lies at distance $0.315$ in \eqref{eq:ksnorm} from the refined
center, measured at 160 temporal harmonics and 60 sine modes rather than the
32 harmonics and 48 modes of Table~\ref{tab:stations}. This candidate comes
from an earlier run of the same FD2, $M=24$ pipeline (held-out error
$6.4\times10^{-4}$), not from the sweep in Table~\ref{tab:stations}, whose
$M=24$ FD2 candidate converges to a different orbit of period $73.12$. The
chaotic reference trajectory is sensitive to floating-point details, so the
two runs' data and fitted models differ slightly. Part of the distance $0.315$
is the 32-harmonic temporal truncation of the candidate, and part is the
error of the surrogate. The certificate establishes
the existence of a true orbit found through the station surrogate. Unlike
Theorem~\ref{thm:ksorbit}, it does not quantify the accuracy of the frozen
surrogate export.

\paragraph{Station implementation.}
The sensor-count study is \path{CAP/Kuramoto_PO/sensor_convergence.py}, with
the table and figure produced by \path{summarize_sensor_convergence.py}. The
reconstruction routines and the irregular-site comparison are in
\path{sensor_experiment.py}. The refined station certificate is produced by
\path{sensor_proof_candidate.py} and \path{run_sensor_refined_proof.sh}.

\FloatBarrier
\paragraph{Certification implementation.}
Candidate exports, certificate data, and proof drivers are provided in
\path{CAP/Mackey-Glass_PO/} and \path{CAP/Kuramoto_PO/} in the
\href{https://github.com/Aloh5/learned-dynamics-cap}{project's GitHub repository}.
The accompanying documentation specifies the MATLAB/INTLAB requirements
and the arithmetic assumptions used for validation.

\section{Discussion and conclusions}
\label{sec:discussion}
The proposed validation method uses certified invariant objects to assess a
learned model's ability to reproduce selected structures of its reference
system. Its mathematical and computational ingredients are established; the
methodological contribution is to use their combination as a model-validation
procedure with an explicit candidate origin, error enclosure, and scope of
interpretation. The two applications demonstrate that this test can succeed
for models of delay and partial differential equations. All physical
parameters remain fixed, and all omitted modes are included. The two
certificates centered at surrogate exports require no true-equation Newton
corrections. Separate refinements use five Newton updates for Mackey--Glass
and three for Neumann KS, with certified radii $1.2\times10^{-12}$ and
$7\times10^{-9}$ around the respective corrected centers. Each refined
ball is contained in its uncorrected enclosure in the same phase slice,
identifying the same orbit. These results quantify refinement accuracy;
the uncorrected enclosures measure the exported candidates' accuracy.

The forecasting and qualitative results remain essential to the assessment.
The cubic Mackey--Glass model supports useful forecasts, long irregular
integrations, and a close approximation to the leading equilibrium
instability, while also possessing additional equilibria. Its attractor
projections and empirical distributions supply numerical evidence beyond
pointwise trajectory agreement. The degree-five refinement model gives a
different outcome: its global rollouts fail, yet it supplies a certifiable
periodic candidate. The structural validation detects a particular success
without erasing the observed failures.

For periodic KS, the forecasting, equilibrium, and attractor comparisons
concern a 14-node spectral semidiscretization. The rigorous orbit result
concerns a separate Neumann problem and cannot certify the periodic-domain
forecasts. Near-machine-precision derivative recovery in the Neumann fit
reflects the matching quadratic library and exact derivative supervision.
Candidate construction and certification test additional properties of that
fit. These distinctions allow the numerical and rigorous results to contribute
to model assessment without treating any single diagnostic as sufficient.

The station experiment of Section~\ref{sec:station-results} removes the exact
derivative supervision. With up to 24 sensors, quadratic surrogates learned from
station time series produce periodic candidates. The error of those candidates
decreases geometrically with the number of sensors, every candidate converges
to a true orbit within two to nine Newton corrections, and one such orbit is
certified after refinement (Theorem~\ref{thm:ksorbit-stations}). At $M=24$ the
remaining error comes mostly from unresolved modes, not from the learned
dynamics. Adding sensors to resolve those modes makes the regression unstable,
because data sampled along the attractor carry almost no information about
the damping of the fine modes. The diagnostic therefore exposes a limitation of
the polynomial surrogate class that held-out derivative errors do not reveal.

\subsection{Limitations and extensions}
These examples establish feasibility of the validation method. They do not
establish superiority over other validation approaches or learned models, or a
population certification rate. A candidate ensemble with specified search
and proof budgets would be needed to measure success fractions. Proof radii
also depend on norms, scaling, conditioning, and the sharpness of estimates;
they cannot be compared directly across these two systems.

Reference-integrator and spatial-resolution studies remain necessary for
stronger forecasting conclusions. Exact-derivative operator recovery does not
address noisy or incomplete observations; Section~\ref{sec:future-work}
outlines a proposed robustness study. Dense libraries grow as $\binom{d+p}{p}$,
and numerical survival supplies no global trapping guarantee. Structural
constraints could improve long-time behavior, while validated continuation
and Floquet analysis could establish parameter dependence and orbital
stability, neither proved here.

A connection to inertial-manifold reduction is also worth investigating
\cite{foias1988,zelik2014}, but a Galerkin fit does not reconstruct unresolved
modes as functions of retained coordinates. Such a reduction would require a
reconstruction map, uniform invariance-defect estimates, and appropriate
attraction or spectral hypotheses. The present certificates establish none
of those global properties.

\subsection{Future work}
\label{sec:future-work}
The station experiment motivates our next study. A complete quadratic in the
station coordinates cannot be both stable and accurate to the proof radius,
since the regression must identify the damping of modes that the data barely
excite. Neural-network surrogates trained on the same station data offer
different trade-offs: architectures can enforce dissipation or stability, and
can represent the unresolved modes as nonlinear functions of resolved ones
rather than as independent variables. We plan to apply the same pipeline to
such surrogates: harmonic balance, frozen export, initial defect, Newton count,
and certification. Other remedies within the polynomial class include learning
a static slaving map for the fine modes from the same data and fitting the
coefficients of a PDE ansatz instead of an ODE. The latter assumes more about
the equation than the present experiments do.

Future experiments should train surrogates on real-world data, including
partially observed and irregularly sampled time series. The present experiments
use synthetic data without added measurement noise and do not establish noise
robustness. A systematic study should add noise to observations before
constructing delay vectors and derivative targets, vary noise levels and
sampling rates, and repeat each experiment with independent noise realizations.
Noise affects both derivative estimates and nonlinear features; ridge
regularization alone does not resolve this errors-in-variables problem.

Weak or integral formulations and ensemble methods provide approaches to
investigate \cite{messengerWeak2021,messenger2021,fasel2022}. Dense ridge,
sparse, and weak-form fits could be compared using separate training and
validation data, clean held-out forecast targets, long-time diagnostics,
and candidate-certification success under equal search and proof budgets.
Repeated experiments across initial conditions, missing-data patterns, and
datasets would quantify variability. Certifying a particular candidate in a
specified equation would not itself establish robustness to observation noise
or uncertainty in physical parameters.

Additional directions include imposing appropriate physical constraints or
symmetries, testing spatial-resolution and reference-integrator convergence,
and measuring certification rates for candidate ensembles. Validated
continuation and Floquet analysis could extend the proofs to parameter
dependence and orbital stability. Applications to measured data would also
require distinguishing observation uncertainty and model error from rigorous
validation of a specified source equation.

\subsection{Computational reproducibility}
The repository's \path{CAP/Mackey-Glass_PO/} and \path{CAP/Kuramoto_PO/}
directories contain candidate data, validation drivers, certificates, and proof
audits; the KS materials also contain fitted-model continuation and finite/tail
scans. The certificates record the norms, radii, periods, and correction counts
reported here. Numerical tests check derivatives, Fourier conversions, source
identity, and certificate consistency, but do not replace rigorous validation.
The Newton-refined Mackey--Glass result is reproduced by
\path{CAP/Mackey-Glass_PO/run_newton_refined_proof.sh}; its certificate,
containment check, and candidate history are stored in the accompanying
\path{outputs_newton_refined/} directory. The Neumann KS refinement is
reproduced by \path{CAP/Kuramoto_PO/run_newton_refined_proof.sh}, with its
own \path{outputs_newton_refined/} artifacts and a documented transfer of
the original full-space operator bounds.
Reproduction requires the documented MATLAB/INTLAB environment and arithmetic
assumptions.

\paragraph{Implementation and code availability.}
For the implementations in this paper, we used the code in the GitHub
repository \url{https://github.com/Aloh5/learned-dynamics-cap}.
The regression and numerical assessment workflows are provided in
\path{notebooks/}, and the candidate data and proof drivers are provided in
\path{CAP/Mackey-Glass_PO/} and \path{CAP/Kuramoto_PO/}. The proof computations
require the documented MATLAB/INTLAB environment and arithmetic assumptions.

\enlargethispage{\baselineskip}
\paragraph{Conclusion.}
Learned polynomial models can be evaluated through forecasts, numerical
dynamical diagnostics, and rigorous comparison of selected invariant
structures. The proposed candidate-based validation method supplies the third
assessment by certifying a nearby exact solution around a frozen surrogate
export. The Mackey--Glass and KS results demonstrate its use alongside the
other diagnostics, yielding a more specific account of what each model
reproduces and where its limitations remain.

\begin{samepage}
\paragraph{Acknowledgment of AI assistance.}
AI tools were used to assist with code and repository organization and with a
brief literature review, including locating relevant references and summarizing
related approaches. This assistance supported the research workflow; the authors
remain responsible for the scientific claims, numerical results, mathematical
proofs, and final presentation.

\end{samepage}

\FloatBarrier
\renewcommand{\bibfont}{\small}

\end{document}